%% file: main.tex
\documentclass[Vancouver,Times1COL]{WileyNJDv5} 

\articletype{Article Type}%

\startpage{1}

\input{ex_shared}
\begin{document}

\title{Space-Time Block Preconditioning for Incompressible\\ Resistive Magnetohydrodynamics}

\author[1]{Federico Danieli}

\author[2]{Ben S. Southworth}

\author[3]{Jacob B. Schroder}

\authormark{Danieli \textsc{et al.}}
\titlemark{Space-Time Block Preconditioning for Incompressible Resistive MHD}

\address[1]{\orgdiv{Mathematical Institute}, \orgname{University of Oxford}, \orgaddress{\state{Oxford}, \country{UK}}}

\address[2]{\orgname{Los Alamos National Laboratory}, \orgaddress{\state{New Mexico}, \country{USA}}}

\address[3]{\orgdiv{Department of Mathematics and Statistics}, \orgname{University of New Mexico}, \orgaddress{\state{New Mexico}, \country{USA}}}

\corres{Corresponding author Jacob B. Schroder, \email{jbschroder@unm.edu}}

\presentaddress{Dept. of Mathematics and Statistics; MSC01 1115; University of New Mexico, Albuquerque, NM 87131}

\fundingInfo{EPSRC Centre for Doctoral Training (F. Danieli) in Industrially
Focused Mathematical Modelling (EP/L015803/1), in collaboration with the Culham
Centre for Fusion Energy. B. Southworth was supported as a Nicholas C.
Metropolis Fellow under the Laboratory Directed Research and Development
program of Los Alamos National Laboratory. J. Schroder was supported by NSF
grant DMS-2110917. }

\abstract[Abstract]{
This work develops an all-at-once space-time preconditioning approach for
resistive magnetohydrodynamics (MHD).
We consider parallel-in-time due to the long time
domains required to capture the physics of interest, as well as the complexity
of the underlying system and thereby computational cost of long-time
integration. To ameliorate this cost by using many processors, we thus develop
a novel approach to solving the whole space-time system that is parallelizable
in both space and time. We develop a space-time block preconditioning for
resistive MHD,
following the space-time block preconditioning concept first introduced by
Danieli \emph{et al.} in 2022 for incompressible flow, where an effective
preconditioner for classic sequential time-stepping is extended to the
space-time setting.
The starting point for our derivation is the continuous Schur complement
preconditioner by Cyr \emph{et al.} in 2021, which we proceed to generalise
in order to produce, to our knowledge, the first space-time block
preconditioning approach for the challenging equations governing
incompressible resistive MHD. 
The numerical results are promising for the model problems of island
coalescence and tearing mode, \comment{which we investigate numerically in a mixed regime, where advective and diffusive forces are both present.}  The overhead computational cost
associated with space-time preconditioning versus sequential time-stepping
is modest and primarily in the range of 2$\times$--5$\times$, which is
low for parallel-in-time schemes in general.  Additionally, the scaling
results for inner (linear) and outer (nonlinear) iterations are flat in the
case of fixed time-step size and only grow very slowly in the case of
time-step refinement.  }

\keywords{ Parallel-in-time integration, block preconditioning, magnetohydrodynamics, finite element method }

\maketitle


\section{Introduction}
Lying at the interface between electromagnetism and fluid dynamics, \emph{Magnetohydrodynamics} (MHD) is responsible for modelling the behaviour of electrically charged fluids subject to electromagnetic interactions.
\comment{In our work, we consider \emph{resistive} plasmas as the target application, due to the relevance this covers for our partner institution CCFE \cite{CCFE}.}
\comment{Due to the complexity and cost associated with setting up large-scale experiments, fusion researchers rely extensively on the results from numerical simulations.}

The system of partial differential equations (PDEs) representing resistive MHD is tightly coupled, non-symmetric, and includes nonlinearities that give rise to turbulent behaviour and spawn a range of phenomena interacting at various scales of space and time. 
\comment{For instance, the governing equations allow for the topology of the magnetic field lines to change, leading to phenomena in which these break or coalesce \cite{introPlasma}.}
These are all characteristics that make an approximate numerical solution particularly challenging. Nonetheless, there is great interest in developing strategies for \comment{computationally solving} these problems effectively, as testified by the amount of research dedicated to the design of stable and efficient schemes and preconditioners for accelerating the solution of the discretised MHD system; see for example \cite{MHDstable,MHDPatrick,MHDPatrick2,MHDPrecon1,MHDPrecon2,MHDPrecon3,MHDPrecon4}, and the work this paper primarily builds on \cite{MHDwavePrecon}.
Given the difficulties already posed by having to solve this problem in a time-stepping framework, perhaps it is not surprising that few papers have considered parallel-in-time (PinT) approaches. To our knowledge, only \cite{baudron1,samaddar1,samaddar2} and more recently \cite{samaddar3}, have tested the application of a PinT algorithm (specifically, Parareal \cite{pararealOriginal}) to MHD problems, showing relatively modest speedups of $8$-$10\times$ using hundreds of processors.  However as noted in these works, since the overall time-to-solution for MHD simulations can be prohibitively long (in light of the complexity of the equations involved), even small speedups can have meaningful effects in reducing total time-to-solution, thus making parallel-in-time approaches particularly appealing.

In this paper, we consider a different approach than Parareal to time-parallelisation, based on the principle of \emph{space-time block preconditioning} applied to the all-at-once space-and-time solution of resistive, incompressible MHD. Recently introduced in \cite{danieli2021spacetime}, space-time block preconditioning leverages block preconditioning principles developed for accelerating implicit linear solves in classical time-stepping procedures, and generalises them to the all-at-once solution of space-time systems. In our work, the starting point for defining our space-time block preconditioner is the \emph{Continuous Schur Complement preconditioner} for resistive MHD introduced in \cite[(3.32)]{MHDwavePrecon}. Noticeably, it relies on the very same pressure convection-diffusion (PCD) approach \cite{PCDoriginal,andyFIT}, which is the foundation of the space-time block preconditioner for incompressible flow derived in \cite{danieli2021spacetime}. As such, when considering its extension to the whole space-time framework, \comment{we can use a similar approach as that developed in \cite{danieli2021spacetime}.} 

\comment{Thus, our main contribution is the presentation of this novel space-time block preconditioning approach for MHD, which is a first for this application. We then proceed to validate the new method numerically for some MHD model problems in a mixed regime where advective and diffusive forces are both present. The convergence results are an important first step toward scalable space-time block preconditioners for MHD and also compare well with the preconditioner's single time-step counterpart.  We leave the regime where advective forces dominate as future work. }

The paper proceeds as follows. In \cref{sec::MHD::problem}, we introduce the resistive incompressible MHD system of partial differential equations (PDEs), together with a brief derivation of the model, a justification for the choice of the specific formulation employed, and details of the chosen space-time discretisation. The proposed preconditioning strategy is presented in \cref{sec::MHD::prec}, where we show how to extend the spatial block preconditioning techniques developed in \cite{MHDwavePrecon} from a single time-step to the whole space-time setting. Finally, results on the effectiveness of the preconditioner from a theoretical and experimental point of view are presented and discussed in \cref{sec::MHD::results}.

\input{sec/problem}
\input{sec/precon}
\input{sec/results}

\section{Conclusions}

In this work, we develop a \comment{novel} all-at-once space-time preconditioner for
resistive MHD.
In formulating our preconditioner, we follow the space-time block preconditioning
principle first outlined in \cite{danieli2021spacetime}. There, an application to 
incompressible flow is considered, where a suitable block preconditioning
strategy for sequential time-stepping (namely, PCD \cite{ElmanPCD}) is extended
to the all-at-once, space-time setting. In this work, we showcase how this same
principle can be applied to more complicated PDEs. In particular, we focus on
applications in resistive MHD, and extend the continuous Schur complement
preconditioner introduced in \cite{MHDwavePrecon} from the sequential time-stepping
to the space-time setting. We thus produce, to our knowledge, the first
space-time block preconditioning approach for this problem (or any problem
aside from Navier Stokes). During our derivation, we consider a number
of simplifications and approximations from a direct block LU decomposition
to arrive at a practical space-time preconditioner, which only requires the
parallel solution of inner time-dependent single-variable equations, rather
than the original time-dependent and nonlinearly coupled sets of equations. 

The numerical results \comment{focus on a mixed regime where advective and diffusive forces are both present} 
and are promising, with the overhead associated with
space-time preconditioning versus sequential time-stepping modest in Table \ref{tab::MHD::PCompSingleStep},
typically in the range of $2\times$--$5\times$.  We note that such overhead is
small for a parallel-in-time scheme \cite{ong-schroder-2020}. Additionally, the
scaling results for inner (linear) and outer (nonlinear) iterations are flat in
the case of fixed $\Delta t$ and only very slowly growing in the case of
refinements in time, providing further evidence for the quality of the
preconditioner proposed. Our experiments here rely on using exact block
inverses for the inner scalar (but still space-time) equations, as this initial work 
focuses on the design of the overall space-time block preconditioning approach,
which is a first for MHD. Future work 
\comment{will focus on the regime where advective forces dominate and also study} 
the development of suitable
iterative schemes for the individual
space-time blocks, especially $F_{\bb{u}}$ and $\tilde{C}_A$, where we will leverage
previous work on nonsymmetric AMG \cite{sivas2020air,AIR,lAIR}.

\bibliographystyle{wileyNJD-Vancouver} 
\bibliography{main}

\end{document}

%% file: ex_shared.tex
\usepackage{lipsum}
\usepackage{graphicx}
\usepackage{epstopdf}

\usepackage{hyphenat}
\usepackage{amssymb, amsmath, amsfonts, amsopn, array, empheq}
\usepackage{mathtools, xfrac, scalerel}
\usepackage[makeroom]{cancel}
\usepackage[export]{adjustbox}

\usepackage{hyperref}
\hypersetup{
  colorlinks,
  frenchlinks=false,
  linkcolor={blue!50!black},
  citecolor={blue!50!black},
  urlcolor=black,
  breaklinks=true
  naturalnames=false,
  hypertexnames=false,
}
\usepackage[all]{hypcap}    
\usepackage{cleveref}
\usepackage{enumitem}
\newlist{problems}{enumerate}{1}
\setlist[problems,1]{
    label={\textbf{Problem~\arabic*}},
    ref=\arabic*, 
    wide, labelindent=\parindent, itemsep=\parskip} 
\creflabelformat{problemsi}{#2\textup{#1}#3} 
\crefname{problemsi}{problem}{problems} 
\Crefname{problemsi}{Problem}{Problems}
\crefname{enumi}{Problem}{Problems}     
\newlist{steps}{enumerate}{1}
\setlist[steps,1]{
    label={\textbf{Step~\arabic*}},
    ref=\arabic*, 
    wide, labelindent=\parindent, itemsep=\parskip} 
\creflabelformat{stepsi}{#2\textup{#1}#3} 
\crefname{stepsi}{step}{steps} 
\Crefname{stepsi}{step}{steps}
\crefname{enumi}{Step}{Steps}     
\crefformat{equation}{(#2#1#3)} 
\crefrangeformat{equation}{(#3#1#4) to (#5#2#6)}
\crefmultiformat{equation}{(#2#1#3)}{ and (#2#1#3)}{, (#2#1#3)}{ and (#2#1#3)}
\crefrangemultiformat{equation}{((#3#1#4)) to ((#5#2#6))}{ and (#3#1#4) to (#5#2#6)}{, (#3#1#4) to (#5#2#6)}{ and (#3#1#4) to (#5#2#6)}
\crefname{section}{Sec.}{Sec.}
\Crefname{section}{Section}{Sections}
\crefname{table}{Tab.}{Tab.}
\Crefname{table}{Table}{Tables}
\crefname{figure}{Fig.}{Fig.}
\Crefname{figure}{Figure}{Figures}
\crefname{algorithm}{Alg.}{Alg.}
\Crefname{algorithm}{Algorithm}{Algorithms}
\crefname{theorem}{Thm.}{Thm.}
\Crefname{theorem}{Theorem}{Theorems}
\crefname{remark}{Rmk.}{Rmk.}
\Crefname{remark}{Remark}{Remarks}

\DeclareMathOperator*{\diag}{diag}

\newcommand{\comment}{\textcolor{black}}

\newcommand{\bb}[1]{{\boldsymbol{#1}}}
\newcommand{\mc}[1]{{\mathcal{#1}}}

\usepackage{multirow}
\usepackage{eqparbox}
\newcommand{\tableIndices}[2]{{
  \begin{array}{@{}r@{}}
    \scriptstyle #2\smash{\eqmakebox[ind]{$\scriptstyle\rightarrow$}} \\[-\jot]  
    \scriptstyle #1\smash{\eqmakebox[ind]{$\scriptstyle\downarrow$}}
  \end{array}}}
\newcommand\TopSp{\rule{0pt}{2.6ex}}         
\newcommand\BotSp{\rule[-0.9ex]{0pt}{0pt}}   

\usepackage{tikz,tikz-3dplot}
\usepackage{pgfplots,pgfplotstable}
\usepgfplotslibrary{patchplots,colorbrewer}
\pgfplotsset{compat=1.15}
\usepackage{arydshln,subcaption}
\usetikzlibrary{math, matrix, fadings, shapes, calc, fit, backgrounds, arrows,arrows.meta, decorations.pathreplacing, decorations.markings, positioning}
\tikzset{cross/.style={cross out, draw=black, minimum size=2*(#1-\pgflinewidth), inner sep=0pt, outer sep=0pt},
cross/.default={1pt}}
\usepackage{listofitems,siunitx,enumerate,tabu}
\pgfdeclarelayer{background}
\pgfdeclarelayer{foreground}
\pgfsetlayers{background,main,foreground}   
\usetikzlibrary{external}
\pgfkeys{
    /tikz/external/mode=list and make
}

%% file: sec/problem.tex
\begin{section}{Problem definition and discretisation}
\label{sec::MHD::problem}

This paper focuses on the incompressible resistive MHD model,
described by the system of PDEs
\begin{subequations}
  \label{eqn::MHD::problem::IMHD+}
  \begin{empheq}[left={\empheqlbrace\,}]{align}
  \frac{\partial \bb{u}(\bb{x},t)}{\partial t} + (\bb{u}(\bb{x},t)\cdot\nabla)\bb{u}(\bb{x},t)- \mu \nabla^2\bb{u}(\bb{x},t)\quad&\notag\\
               + \nabla p(\bb{x},t) + \nabla \comment{\left( A(\bb{x},t) j(\bb{x},t) \right)}                &\,=\,\bb{f}(\bb{x},t),
              \label{eqn::MHD::problem::IMHD+::momentum}\\
                              - \nabla \cdot \bb{u}(\bb{x},t)                        &\,=\,0,               \label{eqn::MHD::problem::IMHD+::incomp}\\
                              j(\bb{x},t) - \frac{1}{\mu_0}\nabla^2A(\bb{x},t)       &\,=\,0,                \label{eqn::MHD::problem::IMHD+::current}\\
  \frac{\partial A(\bb{x},t)}{\partial t} + \bb{u}(\bb{x},t)\cdot\nabla A(\bb{x},t)
     - \frac{\eta}{\mu_0}\nabla^2 A(\bb{x},t)                                        &\,=\,-E(\bb{x},t)    \label{eqn::MHD::problem::IMHD+::magnetic},
  \end{empheq}
\end{subequations}
with space-time domain $\Omega\times\left[T_0,T\right]$ and $\Omega\in\mathbb{R}^2$.
   The flow unknowns are given by the velocity and pressure fields, $\bb{u}(\bb{x},t)$ and $p(\bb{x},t)$, while the electromagnetic part of the system is described by the current intensity $j(\bb{x},t)$ and the vector potential for the magnetic field $A(\bb{x},t)$. The latter is defined as $\bb{B}(\bb{x},t) = \nabla\times(A(\bb{x}\comment{,} t)\hat{\bb{k}})$, where $\bb{B}(\bb{x},t)$ denotes the magnetic field and $\hat{\bb{k}}$ represents the unit vector in the $z$-direction. The given vector functions $\bb{f}(\bb{x},t)$ and $E(\bb{x},t)$ represent a forcing term and an externally-imposed electric field, respectively, and act as right-hand sides for our system. The remaining parameters are given by the flow viscosity $\mu$, the electric resistivity $\eta$, and the magnetic permeability in the void $\mu_0$, all of which we consider constant in $\Omega$ for simplicity. 
   Notice \cref{eqn::MHD::problem::IMHD+} is not the only possible formulation for the incompressible resistive MHD system (see for example \cite[Chap.6]{introPlasma}). We consider this choice for several reasons: (i) tracking the magnetic field via its vector potential allows us to strongly impose its solenoidality (note that $\nabla\cdot\bb{B}(\bb{x},t)=\nabla\cdot\nabla\times(A(\bb{x},t)\hat{\bb{k}})\equiv0$); (ii) introducing the current intensity explicitly as a variable in the system prevents us from having to deal with high-order derivatives in the \emph{Lorentz force} term $\nabla \comment{\left( A(\bb{x},t) j(\bb{x},t) \right)}$ appearing in the momentum equation \cref{eqn::MHD::problem::IMHD+::momentum}; and (iii) our work relies on the analysis conducted in \cite{MHDwavePrecon}, which also considers a vector potential formulation for the incompressible MHD equations, similar to equation \cref{eqn::MHD::problem::IMHD+}.

To solve our target system \cref{eqn::MHD::problem::IMHD+} numerically, we apply backward Euler integration in time and a finite-element (FE) discretisation in space. Resorting to an implicit method in time is standard for MHD to ensure stability at relatively large time steps, as system \cref{eqn::MHD::problem::IMHD+} is characterised by dynamics spanning a wide range of time scales. Moreover, backward Euler is also a common choice, because the coupled system is an index-2 DAE, where standard vanilla higher order diagonally implicit Runge-Kutta schemes theoretically limit to first-order accuracy, so higher-order methods in time require additional care. On a uniform temporal grid with spacing $\Delta t=(T-T_0)/N_t$, backward Euler integration gives rise to the discrete nonlinear recurrence relation
\begin{equation}
  \left\{\begin{array}{rcl}
    \frac{\mc{M}_{\bb{u}}}{\Delta t}\left(\bb{u}^{k}-\bb{u}^{k-1}\right)  + \mu\mc{K}_{\bb{u}}\bb{u}^k + \mc{B}^T\bb{p}^k + \mc{N}_{\bb{u}}(\bb{u}^k,\, \bb{j}^k,\bb{A}^k) & = & \bb{f}^k, \\
    \mc{B}\bb{u}^k & = & \bb{0}, \\
    \mc{M}_j\, \bb{j}^k + \frac{1}{\mu_0}\mc{K}_{jA}\bb{A}^k& = & \bb{h}^k, \\
    \frac{\mc{M}_{A}}{\Delta t}\left(\bb{A}^{k}-\bb{A}^{k-1}\right)  + \frac{\eta}{\mu_0}\mc{K}_{A}\bb{A}^k + \mc{N}_{A}(\bb{u}^k,\bb{A}^k) & = & -\bb{E}^k\\
  \end{array}\right.\textnormal{for }\quad k=1,...,N_t,
  \label{eqn::MHD::problem::recurrence}
\end{equation}
where the discrete vector $\bb{u}^k\in\mathbb{R}^{N_\bb{u}}$ represents the (spatial) approximation of the solution $\bb{u}^k\approx\bb{u}(\cdot,t_k)$ at instant $t_k=T_0+k\Delta t$, and similarly for the other unknowns. Notice that the number of degrees of freedom for each variable, $N_{\bb{u}}$, $N_{p}$, $N_{j}$ and $N_{A}$, can differ in general, as these ultimately depend on the specific FE spaces considered for their discretisations: more details on the choice of spaces are given in \cref{sec::MHD::results}.
In \cref{eqn::MHD::problem::recurrence}, $\mc{M}_{(*)}$ and $\mc{K}_{(*)}$ respectively identify the Galerkin mass and stiffness matrices for variable ${(*)}$ ($\mc{K}_{jA}$ is the \emph{mixed} current-vector potential stiffness matrix), $\mc{B}$ denotes the discrete (negative) divergence operator, while $\bb{f}^k$, $\bb{h}^k$, and $\bb{E}^k$ are the corresponding discrete representations of the right-hand sides for the system, including boundary contributions at time $t_k$. The nonlinear operators $\mc{N}_{\bb{u}}$ and $\mc{N}_{A}$ discretise the advection and Lorentz force terms $(\bb{u}\cdot\nabla)\bb{u}+j\nabla A$ in \cref{eqn::MHD::problem::IMHD+::momentum}, and the advection term $\bb{u}\cdot\nabla A$ in \cref{eqn::MHD::problem::IMHD+::magnetic}, respectively.

\subsection{All-at-once solution}

The classical approach to the solution of \cref{eqn::MHD::problem::recurrence} is sequential time-stepping, where the unknowns at time $t_{k+1}$ are (nonlinearly) solved for, given their value at time $t_k$. When resolving the nonlinearities in \cref{eqn::MHD::problem::recurrence} using Newton's method, this gives rise to a sequence of linear systems involving the (approximate) Jacobian of the operator identified by equation \cref{eqn::MHD::problem::recurrence}. This Jacobian for sequential time-stepping can be defined compactly in the following block form:
\begin{equation}
  \mathcal{J}\left(\left[\begin{array}{c}\bb{u}^k\\\bb{p}^k\\\bb{j}^k\\\bb{A}^k\end{array}\right]\right) \coloneqq \left[\begin{array}{c|c|c|c}
    \mc{F}_{\bb{u}}(\bb{u}^k) & \mc{B}^T & \mc{Z}_j(\bb{A}^k) & \mc{Z}_A(\bb{j}^k)         \\\hline
    \mc{B}                    &          &                    &                            \\\hline
                              &          & \mc{M}_j           & \frac{1}{\mu_0}\mc{K}_{jA} \\\hline
    \mc{Y}(\bb{A}^k)          &          &                    & \mc{F}_{A}(\bb{u}^k)       \\
  \end{array}\right].
\label{eqn::MHD::problem::jacobianDiscrete}
\end{equation}
The discrete operator $\mc{F}_{\bb{u}}(\bb{u})$ is composed of the linear terms $\mc{M}_{\bb{u}}/\Delta t$ and $\mu\mc{K}_{\bb{u}}$, as well as the linearisation of $\mc{N}_{\bb{u}}$ in $\bb{u}$ (the latter representing the discretisation of the operator $(\bb{u}\cdot\nabla)(*)+((*)\cdot\nabla)\bb{u})$). Analogously, $\mc{F}_{A}(\bb{u})$ contains $\mc{M}_{A}/\Delta t$ and $\eta/\mu_0\mc{K}_{A}$, plus the linearisation of $\mc{N}_A$ in $\bb{A}$ (that is, the discrete version of operator $\bb{u}\cdot\nabla(*)$). The blocks $\mc{Z}_j(\bb{A}^k)$ and $\mc{Z}_A(\bb{j}^k)$ represent the linearisations of $\mc{N}_{\bb{u}}$ in $\bb{j}$ and $\bb{A}$: these stem from the Lorentz force term, and discretise $(*)\nabla A$ and $j\nabla(*)$, respectively. Finally, $\mc{Y}(\bb{A}^k)$ discretises $\nabla A\cdot(*)$, and comes from the linearisation of $\mc{N}_A$ in $\bb{u}$. 

In this paper, our goal is to recover the solution of \cref{eqn::MHD::problem::recurrence} at \emph{all} discrete temporal and spatial nodes \emph{at once}. This is achieved by solving all the equations composing \cref{eqn::MHD::problem::recurrence} simultaneously for all time instants $k$. The fact that the nonlinearities present in the recurrence relation need to be solved for does not change, and we still resort to Newton's method in order to do so. Furthermore, the full space-time Jacobian maintains the same block structure as in equation \cref{eqn::MHD::problem::jacobianDiscrete}, that is:
\begin{equation}
  J\left(\left[\begin{array}{c}\bb{u}\\\bb{p}\\\bb{j}\\\bb{A}\end{array}\right]\right) \coloneqq \left[\begin{array}{c|c|c|c}
    F_{\bb{u}}(\bb{u}) & B^T & Z_j(\bb{A}) & Z_A(\bb{j})   \\\hline
    B                  &     &             &               \\\hline
                       &     & M_j         & K_{jA}        \\\hline
    Y(\bb{A})          &     &             & F_{A}(\bb{u}) \\
  \end{array}\right].
\label{eqn::MHD::problem::jacobianSTDiscrete}
\end{equation}
The difference in this case is that each of the blocks represents a discrete space-\emph{time} operator, and similarly the variables in this system represent space-time discretisations, obtained by collating the numerical solutions at the set of discrete time values (so that, for example, $\bb{u}\coloneqq\left[(\bb{u}^1)^T,\dots,(\bb{u}^{N_t})^T\right]^T\in\mathbb{R}^{N_tN_\bb{u}}$). In particular, for backward Euler integration, we have 
\begin{equation}
  F_{\bb{u}}(\bb{u})\coloneqq\left[\begin{array}{ccc}
    \mc{F}_{\bb{u}}(\bb{u}^1)         &      &                             \\
    -\frac{\mc{M}_{\bb{u}}}{\Delta t} &\ddots&                             \\
                                      &\ddots&\mc{F}_{\bb{u}}(\bb{u}^{N_t})
  \end{array}\right]\in\mathbb{R}^{N_\bb{u}N_t \times N_\bb{u}N_t},
  \label{eqn::MHD::problem::STFu} 
\end{equation}
and analogously
\begin{equation}
  F_{A}(\bb{u})\coloneqq\left[\begin{array}{ccc}
    \mc{F}_{A}(\bb{u}^1)         &      &                        \\
    -\frac{\mc{M}_{A}}{\Delta t} &\ddots&                        \\
                                 &\ddots&\mc{F}_{A}(\bb{u}^{N_t})
  \end{array}\right]\in\mathbb{R}^{N_AN_t\times N_AN_t}.
  \label{eqn::MHD::problem::STFA} 
\end{equation}
Notice that both equations \cref{eqn::MHD::problem::STFA,eqn::MHD::problem::STFu} present a block {bi}-diagonal structure, as they approximate a temporal derivative via a {one}-step method. Conversely, the differential operators discretised by the remaining blocks composing equation \cref{eqn::MHD::problem::jacobianSTDiscrete} act on the spatial component only, hence their space-time discretisations present a block \emph{diagonal} structure (\emph{i.e.}, they correspond to a set of spatial operators decoupled in time):
\begin{align}
  Z_j(\bb{A})&\coloneqq \diag_{k=1,\dots,N_t}\left(\mc{Z}_A(\bb{A}^k)\right),& Z_A(\, \bb{j})& \coloneqq \diag_{k=1,\dots,N_t}\left(\mc{Z}_j(\bb{j}^k)\right),
  \label{eqn::MHD::problem::lorentz}\\
  M_j        &\coloneqq \diag_{N_t}\left(\mc{M}_j\right),                    & K_{jA}     & \coloneqq \diag_{N_t}\left(\frac{1}{\mu_0}\mc{K}_{jA}\right),\\
             &\qquad\text{and}\qquad                                         & Y  (\bb{A})& \coloneqq \diag_{k=1,\dots,N_t}\left(\mc{Y}(\bb{A}^k)\right),
\label{eqn::MHD::problem::STother}
\end{align}
where $\diag_{N_t}((*))$ and $\diag_{k=1,\dots,N_t}((*)_k)$ both denote $N_t\times N_t$ block diagonal matrices: for the former, the block diagonal is filled with $N_t$ identical blocks containing $(*)$; for the latter, the diagonal blocks are taken in order from the indexed operators $\{(*)_k\}_{k=1}^{N_t}$.

Solving \cref{eqn::MHD::problem::recurrence} in an all-at-once fashion is typically more computationally expensive than doing so via traditional sequential time-stepping; on the other hand, time-stepping is an inherently sequential procedure, \emph{i.e.} the computational work for each time step relies on the previous step, and thus cannot be parallelised. By pursuing an all-at-once approach, we expose additional parallelism in the time dimension. This paper focuses on the construction an efficient, time-parallel solver for \cref{eqn::MHD::problem::recurrence}, with minimal overhead/additional cost when compared with traditional time-stepping. To do so, we propose solving the space-time Jacobian \cref{eqn::MHD::problem::jacobianSTDiscrete} via GMRES, and to accelerate its convergence by means of a preconditioner, designed so that its application can be effectively parallelised over both the spatial and temporal domain.
\end{section}

%% file: sec/precon.tex
\begin{section}{Space-time block preconditioning for incompressible MHD}
\label{sec::MHD::prec}
In designing a preconditioner for \cref{eqn::MHD::problem::jacobianSTDiscrete}, we exploit its block structure and resort to the \emph{space-time block preconditioning} approach first proposed in \cite{danieli2021spacetime} for incompressible Navier Stokes. This approach consists of identifying a block preconditioner that has shown its efficacy applied to single time-step (spatial) operators \cref{eqn::MHD::problem::jacobianDiscrete}, and extending it to be applicable to its space-time counterpart \cref{eqn::MHD::problem::jacobianSTDiscrete}. 

The single time-step spatial preconditioner we start from was initially developed by Cyr \emph{et al.} in \cite{MHDwavePrecon}. The motivation behind this particular choice stems from the specific design of this preconditioner.  It acts by splitting the solution of the fully coupled problem into two separate sub-problems, one dealing with the coupling between velocity and pressure, the other with velocity and magnetism. Solution of the former sub-problem relies on the PCD preconditioner, which served as a building block for designing its space-time equivalent in \cite{danieli2021spacetime}. As we show in this section, analogous ideas to those described in \cite{MHDwavePrecon} can be exploited in the space-time case as well, and result in a similar sub-problem separation. This allows us to reuse techniques developed in \cite{danieli2021spacetime} for the velocity-pressure MHD equations, and to focus our effort on finding an adequate space-time block preconditioner for the velocity and magnetism equations that arise in MHD.

In \cite{MHDwavePrecon}, the authors start from an LDU factorisation of an approximate Jacobian. Adapted to our formulation \cref{eqn::MHD::problem::IMHD+}, and extended to the whole space-time case, this approximation reads
\begin{equation}
\begin{split}
J &\approx \begin{tikzpicture}[
      baseline=(current  bounding  box.center), 
      Highlight/.style={
          draw,
          densely dotted,
          rounded corners=2pt,
      },
      ampersand replacement=\&,
  ]
  \matrix[matrix of nodes,align=center,inner sep =2pt,nodes in empty cells,left delimiter={[},right delimiter ={]}] at (0,0) (A){ 
  $  F_{\bb{u}} $\&$ B^T                 $\&$ Z_j $\&$ Z_A    $\\
  $  B          $\&$                     $\&$     $\&$        $\\
  $             $\&$                     $\&$ M_j $\&$ K_{jA} $\\
  $  Y          $\&$ YF_{\bb{u}}^{-1}B^T $\&$     $\&$ F_{A}  $\\
  };
  \draw[Highlight](A-4-2.north west) rectangle (A-4-2.south east);
	\end{tikzpicture}\\
  &= \underbrace{\left[\begin{array}{cccc}
    F_{\bb{u}} &   & Z_j & Z_A    \\
               & I &     &        \\
               &   & M_j & K_{jA} \\
    Y          &   &     & F_{A}  \\
  \end{array}\right]}_{\eqqcolon P_{\bb{u}jA}}
  \underbrace{\left[\begin{array}{cccc}
    F_{\bb{u}}^{-1} &   &   &  \\
                    & I &   &  \\
                    &   & I &  \\
                    &   &   & I\\
  \end{array}\right]}_{\eqqcolon F_{\bb{u}I}^{-1}}\underbrace{\left[\begin{array}{cccc}
    F_{\bb{u}} & B^T &   &  \\
    B          &     &   &  \\
               &     & I &  \\
               &     &   & I\\
  \end{array}\right]}_{\eqqcolon P_{\bb{u}p}}.
\end{split}
\label{eqn::MHD::prec::jacobianApprox}
\end{equation}
While this factorisation does not correspond to the exact Jacobian \cref{eqn::MHD::problem::jacobianDiscrete} (the error term in the resulting matrix is outlined in equation \cref{eqn::MHD::prec::jacobianApprox}), it has the advantage of neatly separating the original fully-coupled system into two simpler subsystems, denoted as $P_{\bb{u}jA}$ and $P_{\bb{u}p}$. This effectively isolates the influence of certain variables in each system, making the Schur complements simpler to approximate, as described in the following subsections.

\begin{subsection}{Velocity-pressure coupling}
\label{sec::MHD::prec::VP}
The rightmost matrix $P_{\bb{u}p}$ in \cref{eqn::MHD::prec::jacobianApprox} considers only the coupling between the velocity and pressure variables: applying the inverse of this factor is equivalent to solving an incompressible flow problem, without considering the effects of the magnetic field. This system corresponds to a similar problem\footnote{The sole difference is that in our case $\mc{F}_{\bb{u}}$ additionally contains the extra term
\begin{equation}
  \left[\Delta\mc{W}_{\bb{u}}(\bb{u}(\bb{x}))\right]_{m,n=0}^{N_{\bb{u}}-1} = \int_{\Omega} \left(\left(\bb{\phi}_n(\bb{x})\cdot\nabla\right)\bb{u}(\bb{x})\right)\cdot\bb{\phi}_{m}(\bb{x})\,d\bb{x}
  \label{eqn::MHD::prec::advVdelta}
\end{equation}
stemming from the Newton linearisation. As this term includes an interplay between different components of the velocity variable, it is not feasible to build its counterpart on the scalar pressure field, so the commutation argument cannot be followed exactly. To our knowledge, however, the common approach in the single time-step case is to simply ignore this extra contribution to the definition of $\mc{F}_p$ \cite{CyrPCD,ElmanPCD}. We follow the same approach here: notice that the definition of $\mc{F}_p$ in \cref{eqn::MHD::prec::Fp} matches the one in \cite{danieli2021spacetime}.} as the one considered in \cite{danieli2021spacetime}. As such, we follow the same approach to define the corresponding space-time block preconditioner.
Specifically, we first construct the LU factorisation of $P_{\bb{u}p}$,
\begin{equation}
P_{\bb{u}p} = \underbrace{\left[\begin{array}{cccc}
    I              &   &   &  \\
    BF_\bb{u}^{-1} & I &   &  \\
                   &   & I &  \\
                   &   &   & I\\
\end{array}\right]}_{\eqqcolon L_{\bb{u}p}} \underbrace{\left[\begin{array}{cccc}
    F_{\bb{u}} & B^T &   &   \\
               & S_p &   &   \\
               &     & I &   \\
               &     &   & I \\
\end{array}\right]}_{\eqqcolon U_{\bb{u}p}},
\label{eqn::MHD::prec::PupLU}
\end{equation}
and then substitute its block upper triangular term $U_{\bb{u}p}$ with the approximation
\begin{equation}
	\tilde{U}_{\bb{u}p}\coloneqq\left[\begin{array}{cccc}
    F_{\bb{u}} & B^T         &   &   \\
               & \tilde{S}_p &   &   \\
               &             & I &   \\
               &             &   & I \\
	\end{array}\right],
	\label{eqn::MHD::prec::Utilde}
\end{equation}
where the space-time pressure Schur complement $S_p\coloneqq - B F_{\bb{u}}^{-1} B^T$ is approximated by the operator $\tilde{S}_p\coloneqq-K_pF_p^{-1}M_p\approx S_p$. Here, $F_p=F_p(\bb{u})$ is a space-time reaction-advection-diffusion operator acting on the pressure field, and stems from a commutation argument with its velocity counterpart $F_{\bb{u}}$, namely
\begin{equation}
  F_p(\bb{u})\coloneqq\left[\begin{array}{ccc}
    \mc{F}_{p}(\bb{u}^1)         &      &                        \\
    -\frac{\mc{M}_{p}}{\Delta t} &\ddots&                        \\
                                 &\ddots&\mc{F}_{p}(\bb{u}^{N_t})
  \end{array}\right], \quad\text{with}\quad
  \mc{F}_{p}(\bb{u}^{k})\coloneqq\frac{\mc{M}_p}{\Delta t} + \mc{W}_p(\bb{u}^k) + \mu\mc{K}_p,
  \label{eqn::MHD::prec::Fp}
\end{equation}
where $\mc{M}_p$, $\mc{K}_p$, and $\mc{W}(\bb{u}^k)$ are discretisations of pressure mass, stiffness, and advection operators (the latter with velocity field $\bb{u}^k$).
We refer to \cite{danieli2021spacetime} for additional details on its derivation.
\end{subsection}

\begin{subsection}{Velocity-magnetic field coupling}
\label{sec::MHD::prec::VB}

In contrast, the leftmost matrix $P_{\bb{u}jA}$ in \cref{eqn::MHD::prec::jacobianApprox} considers only the coupling between velocity and magnetic fields, omitting the incompressibility constraint associated with the pressure variable. To efficiently precondition the system involving $P_{\bb{u}jA}$, we once again consider its block LU decomposition:
\begin{equation}
P_{\bb{u}jA} = \underbrace{\left[\begin{array}{cccc}
    I                &   &                               &  \\
                     & I &                               &  \\
                     &   & I                             &  \\
    YF_{\bb{u}}^{-1} &   &-YF_{\bb{u}}^{-1}Z_jM_{j}^{-1} & I\\
\end{array}\right]}_{\eqqcolon L_{\bb{u}jA}} \underbrace{\left[\begin{array}{cccc}
    F_{\bb{u}} &   & Z_j   & Z_A    \\
               & I &       &        \\
               &   & M_{j} & K_{jA} \\
               &   &       & S_A    \\
\end{array}\right]}_{\eqqcolon U_{\bb{u}jA}}.
\label{eqn::MHD::prec::jacobianApprox+LeftLU}
\end{equation}
Here, $S_A$ represents the space-time magnetic field Schur complement
of $P_{\bb{u}jA}$, given by
\begin{equation}
  S_A \coloneqq F_{A} - YF_{\bb{u}}^{-1}\underbrace{(Z_A - Z_jM_j^{-1}K_{jA})}_{\eqqcolon Z},
  \label{eqn::MHD::prec::ASchurApprox1+}
\end{equation}
where we compress the linearisation of the Lorentz force term into
$Z(j(\bb{x},t),A(\bb{x},t)) \coloneqq Z_A(j(\bb{x},t)) - Z_j(A(\bb{x},t))M_j^{-1}K_{jA}
= \diag_{k=1,\dots,N_t}\left(\mc{Z}(j(\bb{x},t_k),A(\bb{x},t_k))\right)$,
with individual blocks given by\footnote{
    Notice that in \cite[(3.7)]{MHDwavePrecon}, the single time-step magnetic Schur complement assumes the form
    \begin{equation}
      \mc{S}_A = \mc{F}_{A} - \mc{Y}\mc{F}_{\bb{u}}^{-1}\mc{Z},
      \label{eqn::MHD::prec::ASchurApprox1singleTS}
    \end{equation}
    where $\mc{Z}$ comes from the linearisation of the Lorentz force. In virtue of our formulation \cref{eqn::MHD::problem::IMHD+}, which explicitly considers the current as an additional variable, this term appears in a different form in \cref{eqn::MHD::prec::ASchurApprox1+}, \emph{i.e.}, it is further split into its two components $\mc{Z}_A$ and $-\mc{Z}_j\mc{M}_j^{-1}\mc{K}_{jA}$. Nonetheless, as it represents the same operator, we can still apply similar considerations used in \cite{MHDwavePrecon} for approximating \cref{eqn::MHD::prec::ASchurApprox1singleTS} here.}
\begin{equation}
  \mc{Z}(j(\bb{x}),A(\bb{x}))\coloneqq \mc{Z}_A(j(\bb{x})) - \mc{Z}_j(A(\bb{x}))\mc{M}_j^{-1}\mc{K}_{jA}.
  \label{eqn::MHD::prec::Zequivsingle}
\end{equation}
The main goal of this section is to develop an approximate inverse of the space-time magnetic field Schur complement \cref{eqn::MHD::prec::ASchurApprox1+} that is relatively easy to apply, with particular consideration for the term $YF_{\bb{u}}^{-1}Z$.

\begin{subsubsection}{An inner wave approximation}
\label{sec::MHD::prec::VB::wave}

We now proceed to derive our approximations. 
Analogous to \cite[Sec. 3.1]{danieli2021spacetime}, $YF_{\bb{u}}^{-1}Z$ is composed of blocks in the form
\begin{equation}
  \mc{Y}_k\mc{F}_{\bb{u},k}^{-1}\left(\prod_{l=i}^{k-1}\left(\frac{\mc{M}_{\bb{u}}}{\Delta t}\mc{F}_{\bb{u},l}^{-1}\right)\right)\mc{Z}_i^T,
  \quad\text{with}\quad\left\{\begin{array}{l} k=1,\dots,N_t\\i=1,\dots,k\end{array}\right.,
  \label{eqn::MHD::prec::FinvBlock}
\end{equation}
where $\mc{Z}_i$, $\mc{Y}_i$ and $\mc{F}_{\bb{u},i}$ are shorthand notation for $\mc{Z}(j(\bb{x},t_i),A(\bb{x},t_i))$, $\mc{Y}(A(\bb{x},t_i))$, and $\mc{F}(\bb{u}(\bb{x},t_i))$, respectively, for some fixed fields $\bb{u}(\bb{x},t)$, $j(\bb{x},t)$, and $A(\bb{x},t)$ (note, subscript $i$ in, \textit{e.g.}, $\mc{Z}_i$ corresponds to an arbitrary \emph{index}, as opposed to the specific subscript $\mc{Z}_j$ corresponding to the variable $j(\bb{x})$ \eqref{eqn::MHD::problem::jacobianDiscrete}).

In \cite{MHDwavePrecon} the authors make use of a small-commutator argument akin to the one often exploited in incompressible fluids (e.g., see \cite{PCDoriginal,andyFIT}), and assume that\footnote{Notice that also in this case the extra term \cref{eqn::MHD::prec::advVdelta} in $\mc{F}_{\bb{u},i}$ does not find an equivalent in $\mc{F}_{A,i}$, but we accept this discrepancy under similar considerations as for \cref{sec::MHD::prec::VP}.}
\begin{equation}
  \mc{M}_A^{-1}\mc{F}_{A,i}\mc{M}_A^{-1}\mc{Y}_k \approx \mc{M}_A^{-1}\mc{Y}_k \mc{M}_{\bb{u}}^{-1} \mc{F}_{\bb{u},i},\qquad \forall i,k=1,\dots,N_t.
	\label{eqn::MHD::prec::smallComm}
\end{equation}
We then recursively apply assumption \cref{eqn::MHD::prec::smallComm} to recover reasonable approximations to the blocks \cref{eqn::MHD::prec::FinvBlock}, and thus to the whole space-time operator \cref{eqn::MHD::prec::ASchurApprox1+}. Left-multiplying \cref{eqn::MHD::prec::smallComm} by $\mc{M}_A\mc{F}_{A,i}^{-1}\mc{M}_A$, and right-multiplying it by $\mc{F}_{\bb{u},i}^{-1}$, we recover
\begin{equation}
  \mc{Y}_k\mc{F}_{\bb{u},i}^{-1} \approx \mc{M}_A\mc{F}_{A,i}^{-1}\mc{Y}_k\mc{M}_\bb{u}^{-1},\qquad \forall i,k=1,\dots,N_t.
	\label{eqn::MHD::prec::commutatorfirstblock}
\end{equation}
Starting from $i=k$, continuing to right-multiply by $\frac{\mc{M}_\bb{u}}{\Delta t}\mc{F}_{\bb{u},k-\ell}^{-1}$, increasing the index $\ell$ from $0$ until we reach $k-\ell=l$, repeatedly using \cref{eqn::MHD::prec::commutatorfirstblock} for each $i=k-\ell$, and finally right-multiplying by $\mc{Z}_l$, gives us the following approximation
\begin{equation}
\begin{split}
  &\mc{Y}_k\mc{F}_{\bb{u},k}^{-1}\frac{\mc{M}_\bb{u}}{\Delta t}\mc{F}_{\bb{u},k-1}^{-1} \approx \mc{M}_A\mc{F}_{A,k}^{-1}\underbrace{\frac{\mc{Y}_k}{\Delta t}\mc{F}_{\bb{u},k-1}^{-1}}_{\approx\frac{\mc{M}_A}{\Delta t}\mc{F}_{A,k-1}^{-1}\mc{Y}_k\mc{M}_\bb{u}^{-1}}\\
  \Longrightarrow\quad& \mc{Y}_k\mc{F}_{\bb{u},k}^{-1}\left(\prod_{\ell=1}^{k-l}\left(\frac{\mc{M}_{\bb{u}}}{\Delta t}\mc{F}_{\bb{u},k-\ell}^{-1}\right)\right)\mc{Z}_l
  \approx
	\mc{M}_A\mc{F}_{A,k}^{-1}\left(\prod_{\ell=1}^{k-l}\left(\frac{\mc{M}_A}{\Delta t}\mc{F}_{A,k-\ell}^{-1}\right)\right)\mc{Y}_{k}\mc{M}_\bb{u}^{-1}\mc{Z}_{l},
\end{split}
\label{eqn::MHD::prec::commutatorallblocks}
\end{equation}
which we can apply to each of the blocks in \cref{eqn::MHD::prec::FinvBlock}. In turn, this identifies a first approximation for the space-time magnetic Schur complement \cref{eqn::MHD::prec::ASchurApprox1+}:
\begin{equation}
\begin{split}
	&YF_\bb{u}^{-1}Z \approx M_AF_A^{-1}YM_\bb{u}^{-1}Z\\
  \Longrightarrow\quad&S_A \approx F_{A} - M_A F_{A}^{-1}Y M_{\bb{u}}^{-1} Z = M_A F_{A}^{-1}\left(F_{A}M_A^{-1}F_{A} - YM_{\bb{u}}^{-1}Z\right),
  \label{eqn::MHD::prec::ASchurApprox1}
\end{split}
\end{equation}
which is a direct space-time equivalent of the one proposed in \cite[(3.24)]{MHDwavePrecon} for the single time-step magnetic Schur complement at instant $k$:
\begin{equation}
\begin{split}
  \mc{S}_{A,k} \approx \mc{M}_A \mc{F}_{A,k}^{-1}\left(\mc{F}_{A,k}\mc{M}_A^{-1}\mc{F}_{A,k} - \mc{Y}_k\mc{M}_{\bb{u}}^{-1}\mc{Z}_k\right),\qquad\forall k=1,\dots,N_t.
  \label{eqn::MHD::prec::ASchurApprox1singetimestep}
\end{split}
\end{equation}

An additional approximation is performed in \cite{MHDwavePrecon} in order to further simplify the term $\mc{Y}_k\mc{M}_{\bb{u}}^{-1}\mc{Z}_k$ in \cref{eqn::MHD::prec::ASchurApprox1singetimestep}.  This term in fact proved difficult to treat there, as it heavily disrupts the sparsity pattern of the first term $\mc{F}_{A,k}\mc{M}_A^{-1}\mc{F}_{A,k}$, rendering its solution via multigrid more challenging.\footnote{In our formulation, an additional complication stems from the presence in $\mc{Z}_k$ of the inverse of the mass matrix $\mc{M}_j$, which would also need to be properly approximated (together with $\mc{M}_{\bb{u}}$).} The justification behind the simplification proposed in \cite{MHDwavePrecon} comes from analysing the continuous problem associated with inverting $P_{\bb{u}jA}$, which can be seamlessly extended to the space-time framework. In particular, solving for the velocity-magnetic field coupling \cref{eqn::MHD::prec::jacobianApprox+LeftLU} corresponds to finding a solution to the PDE
\begin{equation}
  \left\{\begin{array}{rcl}
  \displaystyle\frac{\partial \bb{u}}{\partial t} + (\bb{u}\cdot\nabla)\bb{u} - \mu\nabla^2\bb{u} + \frac{1}{\mu_0}\bb{B}\times(\nabla\times\bb{B})&=& \bb{0} \BotSp\\
  \displaystyle\frac{\partial \bb{B}}{\partial t} -\nabla\times(\bb{u}\times\bb{B}) + \frac{\eta}{\mu_0}\nabla\times(\nabla\times\bb{B})           &=& \bb{0} \TopSp
\end{array}\right. \quad\text{in}\quad\Omega\times(T_0,T],
\label{eqn::MHD::prec::reducedMHD}
\end{equation}
obtained by removing the pressure variable in incompressible visco-resistive MHD \cite[see (1)]{MHDwavePrecon}. Solving for the magnetic Schur complement is the discrete equivalent of using the first equation to express $\bb{u}$ in terms of $\bb{B}$, substituting it in the second equation, and recovering its solution $\bb{B}$. To find an approximation to \cref{eqn::MHD::prec::ASchurApprox1singetimestep}, the authors in \cite{MHDwavePrecon} perform a similar procedure, but at the continuous level, and consider a linearised version of \cref{eqn::MHD::prec::reducedMHD}. Using perturbation theory, they expand the solution
$\bb{u}(\bb{x},t) = \bb{u}_0(\bb{x},t) + \bb{u}_1(\bb{x},t)$
and $\bb{B}(\bb{x},t) = \bb{B}_0(\bb{x},t) + \bb{B}_1(\bb{x},t)$
around a quiet state $\bb{u}_0(\bb{x},t)\equiv\bb{0}$ with a constant background magnetic field $\bb{B}_0(\bb{x},t)\equiv \bb{B}_0$. Positioning themselves in the hyperbolic limit, for which the second-order dissipative terms in \cref{eqn::MHD::prec::reducedMHD} are neglected, and dropping high-order perturbations, they recover
\begin{equation}
  \left\{\begin{array}{rcl}
  \displaystyle\frac{\partial \bb{u}_1}{\partial t} + \frac{1}{\mu_0}\bb{B}_0\times(\nabla\times\bb{B}_1) &=& \bb{0} \BotSp\\
  \displaystyle\frac{\partial \bb{B}_1}{\partial t} -\nabla\times(\bb{u}_1\times\bb{B}_0)                 &=& \bb{0} \TopSp
\end{array}\right. \quad\text{in}\quad\Omega\times(T_0,T].
\end{equation}
To find the solution $\bb{B}_1$, they proceed to take an extra time derivative on the second equation, $
  \displaystyle\tfrac{\partial^2 \bb{B}_1}{\partial t^2} -\nabla\times(\tfrac{\partial \bb{u}_1}{\partial t}\times\bb{B}_0)
                                                        -\nabla\times(\bb{u}_1\times\tfrac{\partial \bb{B}_0}{\partial t}) = 0,$
and substitute back in the first equation, to recover
\begin{equation}
\begin{split}
  &\frac{\partial^2 \bb{B}_1}{\partial t^2} - \nabla\times\left(\frac{1}{\mu_0}(\nabla\times\bb{B}_1)\times\bb{B}_0\times\bb{B}_0\right)=0\\
  \Longrightarrow\quad&\frac{\partial^2 \bb{B}_1}{\partial t^2} - \frac{\|\bb{B}_0\|^2}{\mu_0}\nabla^2\bb{B}_1=0.
\end{split}
	\label{eqn::MHD::prec::continuousMagSchur}
\end{equation}
This equation can be analogously expressed in terms of the vector potential,
\begin{equation}
\frac{\partial^2 A}{\partial t^2} - \frac{\|\bb{B}_0\|^2}{\mu_0}\nabla^2A=0,
\label{eqn::MHD::prec::continuousMagSchurA}
\end{equation}
since the formulation in \cref{eqn::MHD::problem::IMHD+::magnetic} implies a unique correspondence between $A$ and its curl. The magnetic Schur complement \cref{eqn::MHD::prec::ASchurApprox1singetimestep} should then approximate, at least to some extent, the action of the operator in equation \cref{eqn::MHD::prec::continuousMagSchurA}, that is, it should represent a wave equation with propagation speed\footnote{This corresponds to the speed at which \emph{Alfv\'{e}n waves} travel: this is a type of magnetohydrodynamic wave generated in response to a perturbation to the magnetic field induced by the presence of currents in the plasma \cite{alfvenWave}.} $\sqrt{\|\bb{B}_0\|^2/\mu_0}$. In \cite{MHDwavePrecon}, this fact is exploited to justify substituting the term $\mc{Y}\mc{M}_\bb{u}^{-1}\mc{Z}$ in equation \cref{eqn::MHD::prec::ASchurApprox1singetimestep} with a discrete Laplacian operator, opportunely scaled by the wave speed. In the space-time case, this translates to the following approximation of the magnetic Schur complement:
\begin{equation}
  S_A \approx M_AF_A^{-1}\underbrace{\left(F_{A}M_A^{-1}F_{A} + K_A\right)}_{\eqqcolon C_A},
  \label{eqn::MHD::prec::ASchurApprox1}
\end{equation}
where we define
\begin{equation}
  K_A\coloneqq\diag_{k=1,\dots,N_t}\left(\frac{\|\bar{\bb{B}}_0^k\|^2}{\mu_0}\mc{K}_A\right),\qquad\text{with}\qquad\bar{\bb{B}}^k \coloneqq \frac{1}{|\Omega|}\int_{\Omega}\nabla\times(A(\bb{x},t_k)\hat{\bb{k}})\,d\bb{x}.
\end{equation}
That is, $\bar{\bb{B}}_0^k$ is taken as the space average of the magnetic field at instant $t_k$.
The task of explicitly assembling $C_A$ in \cref{eqn::MHD::prec::ASchurApprox1} is complicated by the presence of $\mc{M}_A^{-1}$. As a last simplification, we consider instead its approximation $\tilde{C}_A$, where we pick only the diagonal of $\mc{M}_A$, denoted as $\mc{D}_{\mc{M}_A}$, as a surrogate for the whole mass matrix. Overall, the space-time operator $\tilde{C}_A$ has the following block tri-diagonal structure:
\begin{equation}
  \tilde{C}_A\coloneqq\diag_{k=1,\dots,N_t}\left(\tilde{\mc{C}}_{A,k}\right)
 + \frac{1}{\Delta t}  \diag_{k=1,\dots,N_t-1}\left(\tilde{\mc{C}}_{A,-1,k}, -1\right)
 + \frac{1}{\Delta t^2}\diag_{N_t-2}\left(\tilde{\mc{C}}_{A,-2},-2\right),\\
  \label{eqn::MHD::prec::CAspacetime}
\end{equation}
where the notation $\diag_{k=1,\dots,N-i}((*)_k,-i)$ indicates a $N\times N$ block diagonal matrix whose $i$-th block sub-diagonal is filled by taking in order the indexed operators $\{(*)_k\}_{k=1}^{N-i}$. For each block diagonal, these are given by
\begin{subequations}
\begin{align}
  \tilde{\mc{C}}_{A,k}   \coloneqq&\mc{F}_{A,k}\mc{D}_{\mc{M}_A}^{-1}\mc{F}_{A,k} + \frac{\|\bar{\bb{B}}_0^k\|^2}{\mu_0}\mc{K}_A&\quad k=1,\dots,N_t\label{eqn::MHD::prec::CAs::0}\\
  \tilde{\mc{C}}_{A,-1,k}\coloneqq&\mc{M}_A\mc{D}_{\mc{M}_A}^{-1}\mc{F}_{A,k}+\mc{F}_{A,k+1}\mc{D}_{\mc{M}_A}^{-1}\mc{M}_A      &\quad k=1,\dots,N_t-1\label{eqn::MHD::prec::CAs::1}\\
  \tilde{\mc{C}}_{A,-2}  \coloneqq&\mc{M}_A\mc{D}_{\mc{M}_A}^{-1}\mc{M}_A.\label{eqn::MHD::prec::CAs::2}
\end{align}
\label{eqn::MHD::prec::CAs}
\end{subequations}

This allows us to define the final form for our approximate space-time magnetic Schur complement:
\begin{equation}
  S_A \approx \tilde{S}_A \coloneqq M_AF_A^{-1}\tilde{C}_A \qquad\Longleftrightarrow\qquad\tilde{S}_A^{-1} \coloneqq \tilde{C}_A^{-1} F_A M_A^{-1}.
  \label{eqn::MHD::prec::ASchurApprox}
\end{equation}
We thus proceed to simplify the solution of systems involving equation \cref{eqn::MHD::prec::jacobianApprox+LeftLU}, in complete analogy with \cref{sec::MHD::prec::VP}, and substitute its block upper triangular factor $U_{\bb{u}jA}$ with its approximation
\begin{equation}
\tilde{U}_{\bb{u}jA}\coloneqq\left[\begin{array}{cccc}
    F_{\bb{u}} &   & Z_j   & Z_A         \\
               & I &       &             \\
               &   & M_{j} & K_{jA}      \\
               &   &       & \tilde{S}_A \\
\end{array}\right].
\label{eqn::MHD::prec::jacobianApprox+LeftLUapprox}
\end{equation}
\end{subsubsection}

\begin{subsection}{Definition of the space-time block preconditioner}
\label{sec::MHD::prec::P}

In this section, we initially define our \emph{full} space-time block preconditioner for equation
\cref{eqn::MHD::problem::jacobianDiscrete} in \cref{sec::MHD::prec::P::full}, as it stems from
the approximate block LU decomposition outlined in \cref{sec::MHD::prec::VP,sec::MHD::prec::VB}.
We then further simplify this full block preconditioner in \cref{sec::MHD::prec::P::tri}, to obtain
a more efficient block-triangular preconditioner. Using this additional simplification we have observed
comparable performance (in terms of number of iterations to convergence) to the \emph{full} block
preconditioner, but at a smaller cost per iteration, and hence it constitutes our method of choice for our experiments.

\begin{subsubsection}{Full approximation}\label{sec::MHD::prec::P::full}

The full space-time block preconditioner for equation \cref{eqn::MHD::problem::jacobianDiscrete} incorporates the approximate Schur complements introduced in \cref{sec::MHD::prec::VP,sec::MHD::prec::VB} into the approximate factorisation \cref{eqn::MHD::prec::jacobianApprox} and has the form
\begin{equation}
	P\coloneqq L_{\bb{u}jA}\tilde{U}_{\bb{u}jA} F_{\bb{u}I}^{-1} L_{\bb{u}p}\tilde{U}_{\bb{u}p}
	\quad\Longleftrightarrow\quad P^{-1}\coloneqq\tilde{U}_{\bb{u}p}^{-1} L_{\bb{u}p}^{-1}F_{\bb{u}I}\tilde{U}_{\bb{u}jA}^{-1}L_{\bb{u}jA}^{-1},
	\label{eqn::MHD::prec::preconditioner}
\end{equation}
so that applying its inverse involves the following sequence of operations:
\begin{steps}
  \item\label{step::MHDprec::1} \emph{Invert $L_{\bb{u}jA}$}. From its definition in equation \cref{eqn::MHD::prec::jacobianApprox+LeftLU}, after some modifications we can show that applying its inverse to a generic block vector $[\bb{x}_\bb{u}^T,\bb{x}_p^T,\bb{x}_j^T,\bb{x}_A^T]^T$ is equivalent to computing
  \begin{equation}
    L_{\bb{u}jA}^{-1}\left[\begin{array}{c}\bb{x}_\bb{u}\\\bb{x}_p\\\bb{x}_j\\\bb{x}_A\end{array}\right]
    =\left[\begin{array}{c}\bb{x}_\bb{u}\\\bb{x}_p\\\bb{x}_j\\\bb{x}_A+YF_{\bb{u}}^{-1}(Z_jM_{j}^{-1}\bb{x}_j - \bb{x}_{\bb{u}})
\end{array}\right].
  \end{equation}
  This operation requires solving $N_t$ independent systems involving the current intensity mass matrix $\mc{M}_j$, (which can be trivially parallelised over the various time-steps), and more importantly inverting the velocity space-time matrix $F_\bb{u}$ from equation \cref{eqn::MHD::problem::STFu}.

  \item\label{step::MHDprec::2} \emph{Apply $F_{\bb{u}I}\tilde{U}_{\bb{u}jA}^{-1}$}. The two operators can be combined, rendering the task in this step equivalent to inverting the system
  \begin{equation}
    \bar{U}_{\bb{u}jA} \coloneqq \tilde{U}_{\bb{u}jA}F_{\bb{u}I}^{-1} = \left[\begin{array}{cccc}
    I &   & Z_j    & Z_A         \\
      & I &        &             \\
      &   & M_{j}  & K_{jA}      \\
      &   &        & \tilde{S}_A \\
  \end{array}\right].
  \label{eqn::MHD::prec::barUuja}
  \end{equation}
	This step is the most challenging in the application of the preconditioner due to the inversion of the space-time operator $\tilde{S}_A$, and in particular the factor $\tilde{C}_A$ from the wave approximation in \cref{sec::MHD::prec::VB::wave}. From its definition in equation \cref{eqn::MHD::prec::CAspacetime}, we see that this is a fundamentally different operator from $F_\bb{u}$: while the latter essentially represents the space-time discretisation of a standard parabolic PDE, the former rather corresponds to a space-time operator containing high-order mixed spatial and temporal derivatives. 

  \item\label{step::MHDprec::3} \emph{Invert $L_{\bb{u}p}$}. Once again, the main task of this step lies in the inversion of the space-time velocity matrix $F_\bb{u}$, and is similar in complexity to Step 1.

  \item \emph{Invert $\tilde{U}_{\bb{u}p}$}. This task is equivalent to applying the block upper triangular space-time preconditioner $P_T$ proposed in \cite{danieli2021spacetime}, which involves once again inverting the space-time velocity matrix $F_\bb{u}$, as well as the approximate space-time pressure Schur complement $\tilde{S}_p$ (the latter task is trivially parallelisable, as discussed in detail in \cite{danieli2021spacetime}).
\end{steps}

\end{subsubsection}

\begin{subsubsection}{Block-triangular approximation}
\label{sec::MHD::prec::P::tri}

In this section, we now derive our block-triangular approximation, which forms the preconditioner used in practice. 
In block preconditioning $2\times 2$ systems, convergence of Krylov and fixed-point iterations is fully defined by the approximate Schur complement \cite{southworth2020FP}. As a result, in \cite{southworth2020FP} it is proven that one can expect minimal improvement in convergence by using an approximate block LDU compared with an approximate block-triangular preconditioner, at the expense of an additional linear solve of the (1,1) matrix block. Although the analysis in the $4\times 4$ setting equation \eqref{eqn::MHD::problem::jacobianSTDiscrete} is more complicated, we pursue analogous approximations here to reduce the cost of applying the preconditioner.

The heaviest cost in the application of the preconditioner is associated with the solution of the four space-time systems discussed above: three involving $F_\bb{u}$, and one $\tilde{C}_A$. We reduce this number by dropping some of the block lower triangular terms in \cref{eqn::MHD::prec::preconditioner}: this has the result of rendering the preconditioner closer to a block upper triangular form. In particular, by dropping its $L_{\bb{u}jA}$ factor and using the substitution $\bar{U}_{\bb{u}jA}$, we obtain
\begin{equation}
  \tilde{P} \coloneqq \bar{U}_{\bb{u}jA} L_{\bb{u}p}\tilde{U}_{\bb{u}p} = \left[\begin{array}{cccc}
    F_{\bb{u}} & B^T        & Z_j & Z_A          \\
    B          & \Delta S_p &     &              \\
               &            & M_j & K_{jA}       \\
               &            &     & \tilde{S}_{A}\\
  \end{array}\right],
  \label{eqn::MHD::prec::preconditionerSimplified1}
\end{equation}
where $\Delta S_p \coloneqq \tilde{S}_p - S_p$. Notice that applying this simplification amounts to skipping \cref{step::MHDprec::1}, so that one fewer space-time solve for $F_{\bb{u}}$ and one fewer diagonal solve for $M_j$ are necessary for applying the preconditioner. A similar simplification in the single space-time framework is investigated in \cite{MHDwavePrecon} as well.
In addition, we perform one final approximation, and drop also the remaining block lower triangular factor $L_{\bb{u}p}$ (\emph{i.e.}, skip \cref{step::MHDprec::3}), further reducing by one the total number of space-time solves necessary for $F_{\bb{u}}$. With this, we end up with the block upper triangular space-time preconditioner
\begin{equation}
  P_T \coloneqq \bar{U}_{\bb{u}jA}  \tilde{U}_{\bb{u}p} = \left[\begin{array}{cccc}
    F_{\bb{u}} & B^T         & Z_j & Z_A        \\
               & \tilde{S}_p &     &            \\
               &             & M_j & K_{jA}     \\
               &             &     & \tilde{S}_A\\
  \end{array}\right].
  \label{eqn::MHD::prec::PT}
\end{equation}

\end{subsubsection}
\end{subsection}

\begin{subsection}{Comparison with single time-step block preconditioner}
\label{sec::MHD::precon::STvsTS}

   A key measure in studying PinT methods is the \emph{overhead cost}, which we define as the ratio of total computational cost of the PinT method vs. that of sequential time stepping. This value is always $\geq 1$ and, assuming effective parallelisation in time\footnote{\comment{Here, we assume effective parallelisation in time, because in our case, this largely implies access to standard parallel linear algebra operations and access to fast solvers for the space-time blocks we invert.  For such a fast inversion, one could for example use hypre~\cite{HYPRE} and it's library of fast multigrid solvers for classic diffusion-dominated problems and also newer reduction-based methods for advection-dominated problems~\cite{AIR,lAIR}.}}, the number of processors required for a parallel speedup over sequential time stepping is largely defined by the overhead cost. To estimate this cost, we follow similar steps as in \cite{danieli2021spacetime}, noting that the ratio of \comment{total computational cost} in our case is determined by two factors:
\begin{itemize}
   \item[(i)] The computational cost of applying the space-time preconditioner once
      vs. applying the \comment{single-step} spatial preconditioner $N_t$ times (once at each time step),
  \item[(ii)] The total number of preconditioner applications required for each approach (this can be broken into a comparison of average linear iterations/Jacobian solves and total average nonlinear iterations). 
\end{itemize}
\comment{Combined, items (i) and (ii) give an implementation-independent cost model for each approach, which we can compare.}

   In this section, we discuss \comment{item (i), the computational cost of a single iteration. An estimate for the iteration ratios in item (ii) will be derived from numerical results in 
   \cite{danieli2021spacetime} and \cref{sec::MHD::results::STvsTS}.} To simplify our exposition, we focus our analysis on the block upper-triangular preconditioner $P_T$ \cref{eqn::MHD::prec::PT}, but similar considerations hold for \cref{eqn::MHD::prec::preconditioner} and \cref{eqn::MHD::prec::preconditionerSimplified1} as well.

The single time-step spatial preconditioner equivalent of $P_T$ at the $k$-th time-step is given by
\begin{equation}
  \mc{P}_{T,k} \coloneqq \underbrace{\left[\begin{array}{cccc}
    I &   & \mc{Z}_{j,k} & \mc{Z}_{A,k}        \\
      & I &              &                     \\
      &   & \mc{M}_j     & \mc{K}_{jA}         \\
      &   &              & \tilde{\mc{S}}_{A,k}\\
  \end{array}\right]}_{\eqqcolon \bar{\mc{U}}_{\bb{u}jA,k}}\underbrace{\left[\begin{array}{cccc}
    \mc{F}_{\bb{u},k} & \mc{B}^T             &   &  \\
                      & \tilde{\mc{S}}_{p,k} &   &  \\
                      &                      & I &  \\
                      &                      &   & I\\
  \end{array}\right]}_{\eqqcolon \tilde{\mc{U}}_{\bb{u}p,k}},
  \label{eqn::MHD::precon::PTTS}
\end{equation}
where $\tilde{\mc{S}}_{p,k}$ is the single time-step pressure Schur complement approximation \cite{danieli2021spacetime}, and $\tilde{\mc{S}}_{A,k}$ is the single time-step magnetic field Schur complement approximation from \cite[(3.32)]{MHDwavePrecon}, defined as
\begin{equation}
  \tilde{\mc{S}}_{A,k}\coloneqq \mc{M}_A\mc{F}_{A,k}^{-1}\tilde{\mc{C}}_{A,k},
  \label{eqn::MHD::prec::ASchurApproxTS}
\end{equation}
with $\tilde{\mc{C}}_{A,k}$ as in \cref{eqn::MHD::prec::CAs::0}. The two factors in \cref{eqn::MHD::precon::PTTS} respectively mimic the velocity-magnetic field and velocity-pressure couplings discussed in \cref{sec::MHD::prec::VB,sec::MHD::prec::VP}. In \cite{danieli2021spacetime}, it is shown that for the right-most factor in \cref{eqn::MHD::precon::PTTS}, the potential speedup of the parallel method depends on two ratios.  The first ratio is $C_{\mc{F}_{\bb{u}}}N_t / C_{F_{\bb{u}}}$, where $C_{\mc{F}_{\bb{u}}}$ and $C_{F_{\bb{u}}}$ are the costs of inverting $\mc{F}_{\bb{u}}$ and $F_{\bb{u}}$ respectively, and provides an indication of the overall efficiency at which we can solve the whole velocity space-time system $F_{\bb{u}}$.  
\comment{In the context of inverting this key block, this ratio compares the cost of the sequential in time preconditioner at every time-step versus the cost of one iteration of the proposed space-time preconditioner, and thus provides an indication of the expected parallel efficiency.}
The second ratio is $N_{it}^{ST}/N_{it}^0$, where
$N_{it}^{ST}$ and $N_{it}^0$ are the average number of iterations to convergence required by the solver for the space-time and \comment{sequential} time-stepping systems, respectively, and will be further specified with subscript ``L" and ``NL" for linear (inner) and nonlinear (outer) iterations.  This ratio measures the overhead in terms of the number of iterations to convergence associated with preconditioning in space-time vs sequential time-stepping.

We proceed by breaking down the cost comparison for the leftmost factor in \cref{eqn::MHD::precon::PTTS} following a similar analysis. The total cost associated with the application of \cref{eqn::MHD::precon::PTTS} is then given by summing these two contributions.  Inverting $\bar{\mc{U}}_{\bb{u}jA,k}$ requires solving for each of its diagonal blocks, as well as applying operators $\mc{K}_{jA}$, $\mc{Z}_{A,k}$ and $\mc{Z}_{j,k}$. When considering its space-time counterpart $\bar{U}_{\bb{u}jA}$ \cref{eqn::MHD::prec::barUuja}, a similar set of operations must be performed, except we need to deal with space-time operators. We have that $M_j$, $K_{jA}$, $Z_A$ and $Z_j$ are all block-diagonal, so their application (or that of their inverse, in the case of $M_j$) is trivially parallelisable over the time-steps.  Assuming $N_t$ processors are available, we then expect this procedure to require a computational time comparable to that for the single time-step equivalent.
The main difference lies in the block lower-triangular space-time operator $\tilde{S}_A$ \cref{eqn::MHD::prec::ASchurApprox}. By looking at the single factors composing it, its inverse requires: (i) inverting the space-time mass matrix $M_A$ (which is block diagonal, and hence again trivially parallelisable); (ii) applying the space-time vector potential operator $F_A$ (which is block \emph{bi}-diagonal, and hence requires some overhead communication between the processors, but of marginal impact in terms of total cost, following considerations similar to those discussed for $F_p$ in \cite{danieli2021spacetime}); finally, (iii) inverting the space-time operator $\tilde{C}_A$ \cref{eqn::MHD::prec::CAspacetime}, which represents the biggest obstacle in determining the parallel efficiency of the application of the space-time magnetic Schur complement approximation. We denote with $C_{\tilde{C}_A}$ the cost associated with operation (iii), and with $C_{\tilde{\mc{C}}_A}$ the average cost per time-step associated with inverting its single time-step counterpart, $\tilde{\mc{C}}_{A,k}$ in \cref{eqn::MHD::prec::ASchurApproxTS}. Combining this with the considerations above regarding $\tilde{\mc{U}}_{\bb{u}p,k}$, we have that the computational times necessary to apply $P_T$ once, and $\mc{P}_{T,k}$ once \emph{\comment{for each time-step}} are, respectively, $C_{F_{\bb{u}}} + C_{\tilde{C}_A}$ and $(C_{\mc{F}_{\bb{u}}} + C_{\tilde{\mc{C}}_A})N_t$. Therefore, \comment{we measure} whether the space-time approach is competitive over the time-stepping one \comment{by determining} if
\begin{equation}
  \frac{N_{it,L}^{ST}}{N_{it,L}^0}\left( C_{F_{\bb{u}}} + C_{\tilde{C}_A}\right)<\left(C_{\mc{F}_{\bb{u}}} + C_{\tilde{\mc{C}}_A}\right)N_t,
\end{equation}
$N_{it,L}^{ST}$ and $N_{it,L}^0$ being the number of iterations to convergence required by the solver for the linearised system in the two cases, respectively. Note that in the analysis in this section, we are neglecting the overhead cost associated with assembling the relevant matrices for the Jacobian during each Newton iteration. However, this cost should be comparable for serial and space-time, as the matrices involved in the two cases are the same.

\end{subsection}

\end{subsection}

\end{section}

%% file: sec/results.tex
\begin{section}{Results}
\label{sec::MHD::results}
In this section, we investigate the effectiveness of the preconditioner introduced in \cref{sec::MHD::prec} for some model problems describing the arising and evolution of instabilities in resistive plasma.
In the framework of resistive plasmas, the instabilities described here result in modifications to the topology of the magnetic field lines, giving rise to separation and reconnection phenomena, usually associated with the generation of current sheets \cite[Chap.~9.6]{introPlasma}. These phenomena have been a fascinating subject of study over the last decades, and research on this area has undergone considerable effort. Providing a complete overview of the topic is beyond the scope of this work, and we rather restrict ourselves to giving a general description of the dynamics triggered within the problems considered in \cref{sec::MHD::results::modelProblems}. Nonetheless, we refer the interested reader to \cite[Chap.~9 and Chap.~20]{introPlasma}, as well as to the excellent review on the topic by Dieter Biskamp, particularly \cite[Chap.~3-4]{plasmaReconnection}.

We point out that these results represent a first of its kind investigation of space-time preconditioning for resistive plasma, and is in this sense an initial work. \comment{As such, we study numerically the case for $\mu$ and $\eta$ equal to 1.0, where diffusive and advective forces are both present.  We leave the case of small values of $\mu$ and $\eta$, where advective forces dominate, for future work, as these cases typically require \emph{ad-hoc} stabilisation and/or modifications to preconditioners (sequential and space-time) beyond the current scope. Nonetheless, the experiments conducted here serve as proof of concept to validate the effectiveness of the proposed preconditioner framework.}

\begin{subsection}{Model problems}
\label{sec::MHD::results::modelProblems}
The problems introduced in this section have been often used in the literature in order to test algorithms for MHD: see for example \cite{MHDPrecon2,MHDPrecon3} for \cref{pb::MHD::tearing}, and \cite{MHDwavePrecon,MHDPrecon1,MHDPrecon2} for \cref{pb::MHD::island}, whose examples we follow in our set-up. We proceed to describe these in detail next.

\begin{problems}
	\item\emph{Tearing mode.}\label{pb::MHD::tearing}
      This is one of the classic instabilities in resistive plasmas. It arises from the interaction of magnetic fields with opposite orientations: in the presence of resistivity, \comment{the opposite fields are allowed to \emph{reconnect}, under the influence of diffusion, thus progressively weakening until annihilation occurs.} At this stage, the magnetic field lines running in one direction reconnect with those in the opposite direction, and produce a so-called \emph{$x$-point}, where the lines intersect (as opposed to \emph{$o$-points}, around which the lines rotate). This evolution is dictated by nonlinear interactions, and it occurs on a relatively large time-scale; nonetheless, \comment{the tearing instability eventually \emph{saturates}.} 

	The system configuration for this first test-case is as follows. The flow is initially at quiet state,
	\begin{equation}
		\bar{\bb{u}}^0_{TM}(\bb{x}) = 0,
		\label{eqn::MHD::prob::ueqTM}
	\end{equation}
	and the vector potential is at the \emph{Harris sheet equilibrium} \cite[Chap.~4.1]{plasmaReconnection}, that is
	\begin{equation}
		A^{eq}_{TM}(\bb{x}) = \frac{1}{\lambda}\ln\left(\cosh(\lambda y)\right),
		\label{eqn::MHD::prob::AeqTM}
	\end{equation}
   for a given $\lambda$. Note that $*_{TM}$ refers to conditions for the tearing mode problem. This must be sustained by the external electric field
	\begin{equation}
		E^{eq}_{TM}(\bb{x}) = \frac{\eta}{\mu_0}\nabla^2A^{eq}_{TM}(\bb{x}) = \frac{\lambda\eta}{\mu_0}\cosh(\lambda y)^{-2},
	\end{equation}
	which is imposed as a right-hand side for the vector potential equation in \cref{eqn::MHD::problem::IMHD+}. The action of $A^{eq}_{TM}(\bb{x})$ on the flow (via Lorentz force) is counterbalanced by the pressure term
	\begin{equation}
		p^{eq}_{TM}(\bb{x}) =\frac{1}{2\mu_0}\cosh(\lambda y)^{-2} - \bar{p}^{eq}_{TM},
		\label{eqn::MHD::prob::peqTM}
	\end{equation}
	where the constant $\bar{p}^{eq}_{TM}$ is selected to ensure that the pressure has zero mean at equilibrium.
	In order to excite the tearing mode, we perturb the equilibrium by considering the following initial conditions for the vector potential:
	\begin{equation}
		\bar{A}^0_{TM}(\bb{x}) = A^{eq}_{TM}(\bb{x}) - \varepsilon\cos(\pi y)\cos\left(\frac{2\pi x}{L}\right),
	\end{equation}
	with $L$ being the $x$-length of the domain. The resulting solution is periodic in $x\in[-L,L]$, and shows alternating $x$- and $o$-points at saturation; the ``eye'' around the $o$-point is referred to as a \emph{magnetic island}. The solution also presents axial symmetries, so instead of considering the whole spatial domain $[-L,L]\times[-1/2,1/2]$, we can simulate only its top-right corner $\Omega = [0,L]\times[0,1/2]$, and flip its solution around both the $x$- and $y$-axes. As such, we prescribe symmetric boundary conditions on sides $x=0$, $x=L$, and $y=0$; on the top side $y=1/2$, instead, we prescribe Dirichlet BC for $A$ (fixed at the equilibrium value \cref{eqn::MHD::prob::AeqTM}), and slip conditions for $\bb{u}$. 
The parameters used in our problem are $\lambda=5$, $\varepsilon=10^{-3}$, and $L=3$.

	\item\emph{Island coalescence.}\label{pb::MHD::island}
	This problem describes the evolution of two magnetic islands in a plasma, as they could stem from the tearing mode perturbation introduced in \cref{pb::MHD::tearing}.
   \comment{Here, the two magnetic islands merge into one (that is, they coalesce), generating an electric discharge. The specific dynamics depend on the value of the magnetic resistivity, but generally reconnection occurs in a time-scale much smaller than that for \cref{pb::MHD::tearing}.}

	As in \cref{pb::MHD::tearing}, we start from a quiet state,
	\begin{equation}
		\bar{\bb{u}}^0_{IC}(\bb{x}) = 0.
		\label{eqn::MHD::prob::ueqIC}
	\end{equation}
	with the vector potential given by a modification to the Harris sheet equilibrium in equation \cref{eqn::MHD::prob::AeqTM}. Note that $*_{IC}$ refers to conditions for the island coalescence problem. This is at times referred to as \emph{Fadeev} or \emph{Schmid-Burgk equilibrium}, since they pioneered its study in \cite{islandCoalescence1,islandCoalescence2}, or also as \emph{corrugated sheet pinch equilibrium} \cite[Chap.~4.3]{plasmaReconnection}, and is given by
	\begin{equation}
		A^{eq}_{IC}(\bb{x}) = \frac{1}{2\pi}\ln\left[\cosh(2\pi y) + \beta \cos(2\pi x)\right],
		\label{eqn::MHD::prob::AeqIC}
	\end{equation}
	for a given parameter $\beta$. This too must be sustained by an external electric field, in this case given by
	\begin{equation}
		E^{eq}_{IC}(\bb{x}) = \frac{\eta}{\mu_0}\nabla^2A^{eq}_{IC}(\bb{x}) = \frac{\eta}{\mu_0}\frac{2\pi(1-\beta^2)}{(\cosh(2\pi y) + \beta\cos(2\pi x))^2}.
	\end{equation}
	The flow equilibrium is maintained by a counterbalancing pressure,
	\begin{equation}
		p^{eq}_{IC}(\bb{x}) =\frac{1}{2\mu_0}\frac{1-\beta^2}{\left(\cosh(2\pi y) + \beta\cos(2\pi x)\right)^2} - \bar{p}^{eq}_{IC},
		\label{eqn::MHD::prob::peqIC}
	\end{equation}
	where again the quantity $\bar{p}^{eq}_{IC}$ is chosen to ensure the equilibrium pressure is zero-mean valued.
	This equilibrium is perturbed by considering as initial conditions for the vector potential
	\begin{equation}
		\bar{A}^0_{IC}(\bb{x}) = A^{eq}_{IC}(\bb{x}) + \varepsilon\cos\left(\frac{\pi}{2}y\right)\cos(\pi x).
	\end{equation}
	Similarly to \cref{pb::MHD::tearing}, in this case, the resulting solution is periodic in $x\in[0,2]$ and presents a number of symmetries, which we exploit to simplify our simulation and reduce the domain from $[0,2]\times[-1,1]$ to its top-left corner $\Omega = \Omega_\square = [0,1]^2$. We impose symmetry on the $x=0$, $x=1$ and $y=0$ axes, while on $y=1$, we prescribe Dirichlet BC for $A$ (again fixed to the equilibrium \cref{eqn::MHD::prob::AeqIC}), and slip conditions for $\bb{u}$, rendering the flow fully enclosed also in this case. The parameters considered for this problem are $\beta=0.2$ and $\varepsilon=10^{-3}$.

\end{problems}

In our experiments, the choice of finite-element approximation for the discretisation of the functional spaces 
falls on the Taylor-Hood $P_3-P_2$ pair for velocity and pressure \cite[Chap.~17.4]{quarterNM}, while we consider $P_1$ for both the vector potential and the current intensity variables. This stems from considerations on \emph{Finite Element Exterior Calculus} \cite{FEEC,FEEC2}, and from the necessity of ensuring that the two forcing terms in the momentum equation, $\nabla p$ and $j\nabla A$, can balance each other out at equilibrium.

Regarding the parameters for the problems, we put ourselves in a somewhat favourable position, by considering only $\mu=\eta=\mu_0=1$. While this choice has the effect of rendering some of the dynamics described in \cref{pb::MHD::tearing,pb::MHD::island} less evident, it allows us to neglect matters regarding stability under advection-dominated regimes, which is beyond the scope of this work, and rather focus on providing a general indication of the space-time preconditioner performance.

The actual code used for the numerical simulations is an expansion of the one developed for \cite{danieli2021spacetime}: as such, it still relies on MFEM, PETSc, and hypre \cite{MFEM,PETSC,HYPRE} for the linear solvers necessary to tackle the linearised system \cref{eqn::MHD::problem::jacobianSTDiscrete} at each Newton iteration. As an initial guess for the discrete space-time velocity, pressure, and vector potential, we pick the projection of the equilibrium conditions, \cref{eqn::MHD::prob::ueqTM}, \cref{eqn::MHD::prob::peqTM}, \cref{eqn::MHD::prob::AeqTM} and \cref{eqn::MHD::prob::ueqIC}, \cref{eqn::MHD::prob::peqIC}, \cref{eqn::MHD::prob::AeqIC}, onto their respective discrete finite element spaces;
the additional variable $j$ is instead initialised as the solution to its associated discrete system in \cref{eqn::MHD::problem::recurrence}
\begin{equation}
	\mc{M}_j \, \bb{j}^k = -\mc{K}_{jA}\bb{A}^k + \bb{h}^k,\qquad\forall k=1,\dots,N_t,
\end{equation}
where $\bb{A}^k$ is the initial guess on the vector potential at instant $k$, so that the initial residual associated with this variable is effectively $0$.

\end{subsection}

\begin{subsection}{Performance of preconditioner}
This section is dedicated to measuring the effectiveness of the proposed space-time preconditioner for accelerating the iterative solution of the linearised incompressible resistive MHD system \cref{eqn::MHD::problem::jacobianSTDiscrete}. In order to evaluate the preconditioner performance, we focus on tracking the number of iterations necessary to achieve convergence, for both the outer Newton solver and the inner GMRES solver. This parameter represents a more indicative measure of the behaviour of the preconditioner than, \emph{e.g.}, time-to-solution, as the latter strongly depends on implementation choices. This becomes even more significant given the possible challenges represented by \cref{step::MHDprec::2}, and particularly by the inversion of $\tilde{C}_A$ in \cref{eqn::MHD::prec::CAspacetime}, whose effective time parallelisation we do not investigate further in this work.

For similar reasons, in this study we consider only the use of \emph{exact} solvers for inverting the relevant block operators appearing in $P_T$. That is to say, we employ sequential time-stepping for the solution of the space-time operators $F_{\bb{u}}$ and $\tilde{C}_A$, and use LU factorisation for inverting the single blocks in their block diagonals $\mc{F}_{\bb{u},k}$ and $\tilde{\mc{C}}_{A,k}$ (see their definitions in \cref{eqn::MHD::prec::FinvBlock} and \cref{eqn::MHD::prec::CAspacetime}), as well as the vector potential mass matrix $\mc{M}_A$ (appearing within $M_A$ in \cref{eqn::MHD::prec::ASchurApprox}), the current intensity mass matrix $\mc{M}_j$ (present in $M_j$ in \cref{eqn::MHD::prec::jacobianApprox+LeftLUapprox}), and the pressure mass and stiffness operators 
apprearing in equation \cref{eqn::MHD::prec::Utilde}.

We also only present results for the \emph{simplified} block upper-triangular preconditioner $P_T$ in \cref{eqn::MHD::prec::PT}, since early experiments showed negligible difference in convergence behaviour with respect to using its \emph{full} version $P$ from \cref{eqn::MHD::prec::preconditioner}.  This makes $P_T$ a much more appealing choice, given how its associated computational cost is significantly reduced, as per the discussion in \cref{sec::MHD::prec::P}.

\begin{table}[t!] 
\caption{\label{tab::MHD::refinementRes} Iterations to convergence for Newton's method applied for the solution of \cref{pb::MHD::tearing,pb::MHD::island}, for various refinement levels of the spatial and temporal grids (the temporal domain is fixed at $T=1$). The outer solver achieves convergence with a tolerance of $10^{-10}$. The inner GMRES solver is right-preconditioned with $P_T$ \cref{eqn::MHD::prec::PT}, and converges with relative tolerance $10^{-2}$ or absolute tolerance of $10^{-14}$. To invert the relevant operators appearing in $P_T$, we use exact solvers (time-stepping and LU). Inside each column, the number of Newton iterations to convergence is reported on the left, the average number of inner GMRES iterations per outer Newton iteration appears on the right (in parentheses).}
\centering
\begin{tabular}{c|c|cc:cc:cc:cc:cc:cc} 
	Pb &$\tableIndices{\Delta x}{\Delta t}$&\multicolumn{2}{c:}{$2^{-2}$} &\multicolumn{2}{c:}{$2^{-3}$} &\multicolumn{2}{c:}{$2^{-4}$} & \multicolumn{2}{c:}{$2^{-5}$} & \multicolumn{2}{c:}{$2^{-6}$} & \multicolumn{2}{c }{$2^{-7}$} \BotSp\\\hline
  \multirow{6}{*}{1}&$  2^{-2}$&$ 5$&$(7.80)  $&$ 5$&$(9.40)   $&$ 5$&$(11.60)  $&$ 5$&$(13.20)  $&$ 5$&$(13.60)  $&$ 5$&$(16.40) $\TopSp\\	
                    &$  2^{-3}$&$ 5$&$(8.40)  $&$ 5$&$(10.20)  $&$ 5$&$(13.00)  $&$ 5$&$(14.80)  $&$ 5$&$(15.00)  $&$ 5$&$(16.60) $\\		
                    &$  2^{-4}$&$ 5$&$(8.00)  $&$ 5$&$(10.40)  $&$ 5$&$(12.20)  $&$ 5$&$(13.80)  $&$ 5$&$(15.40)  $&$ 5$&$(17.20) $\\		
                    &$  2^{-5}$&$ 5$&$(7.40)  $&$ 5$&$(9.40)   $&$ 5$&$(12.00)  $&$ 5$&$(13.40)  $&$ 5$&$(14.00)  $&$ 5$&$(18.40) $\\		
                    &$  2^{-6}$&$ 5$&$(7.60)  $&$ 5$&$(9.40)   $&$ 5$&$(12.00)  $&$ 5$&$(13.40)  $&$ 5$&$(13.80)  $&$ 5$&$(17.00) $\\		
                    &$  2^{-7}$&$ 5$&$(7.60)  $&$ 5$&$(9.80)   $&$ 5$&$(11.80)  $&$ 5$&$(13.20)  $&$ 5$&$(14.40)  $&$ 5$&$(16.80) $\BotSp\\
  \hline
  \multirow{6}{*}{2}&$  2^{-2}$&$ 4$&$(8.75)  $&$ 5$&$(11.00)  $&$ 5$&$(11.80)  $&$ 5$&$(14.20)  $&$ 4$&$(15.75)  $&$ 4$&$(17.00) $\TopSp\\
                    &$  2^{-3}$&$ 5$&$(9.20)  $&$ 4$&$(11.50)  $&$ 4$&$(13.00)  $&$ 4$&$(14.75)  $&$ 4$&$(16.75)  $&$ 4$&$(18.75) $\\
                    &$  2^{-4}$&$ 4$&$(9.75)  $&$ 4$&$(11.25)  $&$ 4$&$(13.50)  $&$ 4$&$(15.50)  $&$ 4$&$(17.25)  $&$ 4$&$(19.25) $\\			
                    &$  2^{-5}$&$ 4$&$(8.00)  $&$ 4$&$(11.00)  $&$ 4$&$(13.75)  $&$ 4$&$(15.25)  $&$ 3$&$(15.67)  $&$ 3$&$(20.33) $\\			
                    &$  2^{-6}$&$ 3$&$(8.00)  $&$ 3$&$(10.00)  $&$ 3$&$(12.33)  $&$ 3$&$(14.00)  $&$ 3$&$(15.67)  $&$ 3$&$(17.33) $\\			
                    &$  2^{-7}$&$ 3$&$(7.67)  $&$ 3$&$(10.67)  $&$ 3$&$(12.33)  $&$ 3$&$(14.00)  $&$ 3$&$(16.00)  $&$ 3$&$(17.67) $		
\end{tabular}
\end{table}

In \cref{tab::MHD::refinementRes} we report convergence results for both inner and outer solvers, collected for various choices of spatial and temporal grid sizes $\Delta x$ and $\Delta t$, on a fixed temporal domain of size $T=1$. There, we see how the number of Newton iterations necessary for convergence remains largely unchanged, regardless of the mesh chosen. The number of GMRES iterations seems to be increasing very slowly with temporal refinement,
where increasing the number of temporal nodes by a factor of $2^5$ results in slightly more than $\sim2$ times as many iterations). 
\comment{For sequential time-stepping, reducing the time-step size can often improve the preconditioning performance, as the previous time-step solution becomes a better initial guess for the next time-step; however, here we are in a different context, considering a global all-at-once space-time preconditioner.  In this context, a modest increase in iteration counts, as the mesh is refined, is not surprising, because there is no improved initial guess, as in the sequential time-stepping case.}
To \comment{further} put these numbers into perspective, notice that for these experiments the spatial unknowns for the most refined discretisation of \cref{pb::MHD::tearing} amount to
$N_{\bb{u}}=222722$, $N_p=49665$, $N_j=N_A=12545$,
while those for \cref{pb::MHD::island} reach $N_{\bb{u}}=296450$, $N_p=66049$, $N_j=N_A=16641$. Considering that the temporal grids count up to $N_t=128$ points, this results in space-time systems of total size
$(N_\bb{u}+N_p+N_j+N_A)\cdot N_t\approx 38\cdot10^6$
   and $\approx52\cdot10^6$ for the two problems, respectively.

The temporal domain $T=1$ considered for the experiments in \cref{tab::MHD::refinementRes} is however too small to allow for the relevant nonlinear behaviours described in \cref{sec::MHD::results::modelProblems} to manifest themselves into the solution. 
To \comment{understand} how the preconditioner performs when these behaviours become evident, we consider another set of experiments, where we fix the time-step and progressively increase the size of the time domain. The corresponding convergence results are reported in \cref{tab::MHD::fixdtRes}. Overall, they draw a similar picture as those in \cref{tab::MHD::refinementRes}, but with the number of inner iterations showing no growth in this case. Thus, the preconditioner is mildly sensitive in this case to the size of $\Delta t$, but not the number of time-steps.

\begin{table}[t!] 
   \caption{\label{tab::MHD::fixdtRes} Convergence results for Newton iterations applied for the solution of \cref{pb::MHD::tearing,pb::MHD::island}, for various refinement levels of the spatial grid. We fix $\Delta t=0.5$, and progressively enlarge the temporal domain to accommodate for more processors. The solvers used are as in \cref{tab::MHD::refinementRes}. Inside each column, the number of Newton iterations to convergence is reported on the left, the average number of inner GMRES iterations per outer Newton iteration appears on the right (in \comment{parentheses}).}
\centering
\begin{tabular}{c|c|cc:cc:cc:cc:cc:cc}
	Pb &$\tableIndices{\Delta x}{T}$&\multicolumn{2}{c:}{$2^{2}$} &\multicolumn{2}{c:}{$2^{3}$} & \multicolumn{2}{c:}{$2^{4}$} & \multicolumn{2}{c:}{$2^{5}$} & \multicolumn{2}{c:}{$2^{6}$} & \multicolumn{2}{c }{$2^{7}$} \BotSp\\\hline
  \multirow{6}{*}{1}&$     2^{-2}$&$ 5$&$(5.80) $&$ 5$&$(6.20) $&$ 5$&$(6.20) $&$ 5$&$(6.60) $&$ 5$&$(6.40) $&$ 5$&$(6.40) $\TopSp\\		
                    &$     2^{-3}$&$ 5$&$(6.80) $&$ 5$&$(6.80) $&$ 5$&$(6.80) $&$ 5$&$(6.80) $&$ 5$&$(6.80) $&$ 5$&$(6.80) $\\		
                    &$     2^{-4}$&$ 5$&$(6.60) $&$ 5$&$(6.60) $&$ 5$&$(6.60) $&$ 5$&$(6.60) $&$ 5$&$(6.60) $&$ 5$&$(6.80) $\\		
                    &$     2^{-5}$&$ 4$&$(6.00) $&$ 5$&$(6.20) $&$ 4$&$(6.25) $&$ 4$&$(6.25) $&$ 5$&$(6.40) $&$ 5$&$(6.40) $\\		
                    &$     2^{-6}$&$ 4$&$(6.00) $&$ 4$&$(6.00) $&$ 4$&$(6.00) $&$ 4$&$(6.00) $&$ 4$&$(6.25) $&$ 4$&$(6.25) $\\		
                    &$     2^{-7}$&$ 4$&$(5.75) $&$ 4$&$(6.00) $&$ 4$&$(6.00) $&$ 4$&$(6.00) $&$ 4$&$(6.00) $&$ 4$&$(6.00) $\BotSp\\\hline		
  \multirow{6}{*}{2}
                    &$     2^{-2}$&$ 5$&$(6.40) $&$ 5$&$(6.40) $&$ 5$&$(6.40) $&$ 5$&$(6.40) $&$ 5$&$(6.40) $&$ 5$&$(6.40) $\TopSp\\
                    &$     2^{-3}$&$ 4$&$(6.75) $&$ 5$&$(7.00) $&$ 5$&$(7.00) $&$ 5$&$(7.00) $&$ 5$&$(7.00) $&$ 5$&$(7.00) $\\
                    &$     2^{-4}$&$ 4$&$(6.25) $&$ 4$&$(6.25) $&$ 4$&$(6.50) $&$ 4$&$(6.75) $&$ 5$&$(6.80) $&$ 5$&$(6.80) $\\
                    &$     2^{-5}$&$ 4$&$(6.50) $&$ 4$&$(6.50) $&$ 4$&$(6.75) $&$ 4$&$(6.50) $&$ 4$&$(6.50) $&$ 4$&$(6.50) $\\
                    &$     2^{-6}$&$ 3$&$(6.33) $&$ 4$&$(6.75) $&$ 4$&$(6.75) $&$ 4$&$(6.75) $&$ 4$&$(6.75) $&$ 4$&$(6.75) $\\
                    &$     2^{-7}$&$ 3$&$(6.33) $&$ 3$&$(6.33) $&$ 3$&$(6.33) $&$ 3$&$(6.33) $&$ 4$&$(6.75) $&$ 3$&$(6.33) $\\
\end{tabular}
\end{table}

\begin{subsection}{Comparison with sequential time-stepping}
\label{sec::MHD::results::STvsTS}
In this section, we provide an indication of how the number of inner and outer iterations to convergence for the space-time preconditioner $P_T$, ($N_{it,L}^{ST}$ and $N_{it,NL}^{ST}$), compare with those of its single time-step counterpart $\mc{P}_{T,k}$, ($N_{it,L}^0$ and $N_{it,NL}^0$), for the test cases introduced in \cref{sec::MHD::results::modelProblems}.

In order to conduct a balanced comparison between the single time-step and the whole space-time approaches, we adjust the absolute tolerance of the linear solver for the former by scaling it by a factor $\sqrt{N_t}$.  For the whole space-time approach, we do the same for the \emph{outer} solver, but leave the absolute tolerance of the inner solver untouched (this might slightly favour the time-stepping approach, but convergence of GMRES occurred due to the \emph{relative} tolerance being met in almost all of cases). Lastly when time-stepping, we reuse the solution value from the previous time-step as the initial guess for the Newton solver at the following time-step.

\begin{table}[b!] 
   \caption{\label{tab::MHD::PCompSingleStep} Ratios $N_{it,NL}^{ST}/N_{it,NL}^0$ and $N_{it,L}^{ST}/N_{it,L}^0$ \comment{(in parentheses)}, where $N_{it,NL}^{ST}$ and $N_{it,L}^{ST}$ are the total number of outer and inner iterations to convergence for the space-time approach (resulting from the same experiments used to fill \cref{tab::MHD::refinementRes}), while $N_{it,NL}^0$ and $N_{it,L}^0$ represent the average number of outer and inner iterations to convergence per time-step, using sequential time-stepping. For the latter, the inner solver for each time-step is GMRES, right-preconditioned with the single time-step counterpart of the block upper-triangular preconditioner, \cref{eqn::MHD::precon::PTTS}. The outer solver for time-stepping reaches convergence with absolute tolerance $10^{-10}/\sqrt{N_t}$, while the inner solver has the same tolerances used for \cref{tab::MHD::refinementRes}. For some of the most extreme parameters combinations, we noticed that the solution from time-stepping ``freezes'' from a certain (reasonably late) instant onwards: that is, the residual after the Newton update becomes smaller than the prescribed tolerance, thus preventing any other update from being carried out (notice we do not implement mesh adaptivity). To account for this, we compute the averages $N_{it,NL}^0$ and $N_{it,L}^0$ by considering the \emph{effective} time-steps only.}
\centering
\begin{tabular}{c|c|cc:cc:cc:cc:cc:cc}
	Pb &$\tableIndices{\Delta x}{\Delta t}$&\multicolumn{2}{c:}{$2^{-2}$} &\multicolumn{2}{c:}{$2^{-3}$} &\multicolumn{2}{c:}{$2^{-4}$} & \multicolumn{2}{c:}{$2^{-5}$} & \multicolumn{2}{c:}{$2^{-6}$} & \multicolumn{2}{c }{$2^{-7}$} \BotSp\\\hline
  \multirow{6}{*}{1}&$ 2^{-2}$&$ 1.05$&$(1.05) $&$ 1.14$&$(1.18) $&$ 1.31$&$(1.45) $&$ 1.38$&$(1.54) $&$ 1.63$&$(1.97) $&$ 1.66$&$(2.41)$\TopSp\\
                    &$ 2^{-3}$&$ 1.11$&$(1.17) $&$ 1.18$&$(1.29) $&$ 1.27$&$(1.54) $&$ 1.45$&$(1.80) $&$ 1.54$&$(2.15) $&$ 1.65$&$(2.40)$\\
                    &$ 2^{-4}$&$ 1.11$&$(1.21) $&$ 1.21$&$(1.40) $&$ 1.40$&$(1.65) $&$ 1.47$&$(1.82) $&$ 1.50$&$(2.35) $&$ 1.65$&$(2.73)$\\
                    &$ 2^{-5}$&$ 1.11$&$(1.13) $&$ 1.25$&$(1.38) $&$ 1.45$&$(1.79) $&$ 1.55$&$(1.93) $&$ 1.66$&$(2.51) $&$ 1.86$&$(3.50)$\\
                    &$ 2^{-6}$&$ 1.11$&$(1.21) $&$ 1.33$&$(1.50) $&$ 1.48$&$(1.89) $&$ 1.60$&$(2.06) $&$ 1.72$&$(2.71) $&$ 1.91$&$(3.52)$\\
                    &$ 2^{-7}$&$ 1.25$&$(1.42) $&$ 1.29$&$(1.56) $&$ 1.51$&$(1.94) $&$ 1.68$&$(2.21) $&$ 1.77$&$(3.06) $&$ 2.01$&$(3.91)$\BotSp\\
  \hline
  \multirow{6}{*}{2}&$ 2^{-2}$&$ 1.07$&$(1.12) $&$ 1.54$&$(1.73) $&$ 1.57$&$(1.78) $&$ 1.72$&$(2.25) $&$ 1.48$&$(2.13) $&$ 1.59$&$(2.44)$\TopSp\\
                    &$ 2^{-3}$&$ 1.54$&$(1.79) $&$ 1.39$&$(1.76) $&$ 1.60$&$(2.05) $&$ 1.60$&$(2.41) $&$ 1.70$&$(3.25) $&$ 1.79$&$(3.70)$\\
                    &$ 2^{-4}$&$ 1.23$&$(1.53) $&$ 1.45$&$(1.90) $&$ 1.60$&$(2.23) $&$ 1.78$&$(2.68) $&$ 1.95$&$(4.36) $&$ 2.18$&$(5.26)$\\
                    &$ 2^{-5}$&$ 1.45$&$(1.60) $&$ 1.45$&$(1.96) $&$ 1.78$&$(2.67) $&$ 2.03$&$(3.19) $&$ 1.60$&$(3.40) $&$ 1.90$&$(5.36)$\\
                    &$ 2^{-6}$&$ 1.09$&$(1.33) $&$ 1.26$&$(1.63) $&$ 1.55$&$(2.20) $&$ 1.57$&$(2.38) $&$ 1.79$&$(4.21) $&$ 1.98$&$(5.20)$\\
                    &$ 2^{-7}$&$ 1.20$&$(1.39) $&$ 1.41$&$(2.00) $&$ 1.60$&$(2.38) $&$ 1.75$&$(2.83) $&$ 1.88$&$(5.01) $&$ 2.04$&$(5.99)$
\end{tabular}
\end{table}

Given the discussion in \cref{sec::MHD::precon::STvsTS}, 
   it is clear that both ratios $N_{it,L}^{ST}/N_{it,L}^0$ and $N_{it,NL}^{ST}/N_{it,NL}^0$ must be kept as close as possible to $1$ for the space-time preconditioning approach to represent an appealing alternative to sequential time-stepping. \comment{For instance, a value of 1.0 would represent no overhead from the space-time approach, while small values (e.g., $< 10$) will represent a modest overhead that compares favorably to other parallel-in-time schemes~\cite{ong-schroder-2020}.} The results in \cref{tab::MHD::PCompSingleStep} show that this condition is satisfied reasonably well: the ratio of Newton iterations hits a maximum of $\sim2$, while that of GMRES iterations sits between $1$ and $3$ for most cases. As $\Delta t\to 0$, this ratio tends to increase, reaching almost $6$ in certain instances. This is not surprising, as in this case the small time step sizes improve the quality of the initial guess for Newton's method during sequential time-stepping.  

Overall, the results in this section showcase the potential effectiveness of our space-time block preconditioning approach in this more complex MHD case. At least when employing \emph{exact} solvers, the space-time preconditioner requires limited overhead in terms of number of iterations to convergence with respect to its single time-step equivalent. This allows significant room for parallel efficiency, so long as effective PinT methods are developed to solve the single-variable space-time systems appearing within our preconditioner, namely involving operators $F_{\bb{u}}$ and $\tilde{C}_A$. Success has already been shown for solving such single-variable space-time systems in parallel, \emph{e.g.}, \cite{sivas2020air,vandewalleMG,ganderMG}.

\end{subsection}

\end{subsection}


\end{section}